\documentclass[12pt]{amsart}

\usepackage{amssymb,amscd,amsthm,extarrows,verbatim,amsmath,color,fancyhdr,mathrsfs}
\usepackage{fancyvrb}
\usepackage{authblk}
\usepackage{multirow}
\usepackage{amsaddr}
\usepackage{graphicx}
\usepackage{algorithm}
\usepackage{algorithmic}
\usepackage{hyperref}
\usepackage{turnstile}
\usepackage{arydshln}
\usepackage{enumitem}
\usepackage{booktabs}
\usepackage{tikz-cd}
\usetikzlibrary{arrows}
\usepackage{xcolor}
\usepackage{pagecolor}
\newcommand{\jdt}[1]{\color{magenta}{\textbf{[JDT: #1]}}\normalcolor}

\newcommand{\rwrd}[1]{\color{green}{\textbf{[RWRD: #1]}}\normalcolor}

\newcommand{\gsc}[1]{\color{white}{\textbf{[GSC: #1]}}\normalcolor}

\usepackage[letterpaper, left=2.5cm, right=2.5cm, top=2.5cm,
bottom=2.5cm,dvips]{geometry}

\def\arrvline{\hfil\kern\arraycolsep\vline\kern-\arraycolsep\hfilneg}

\newtheorem{lemma}{Lemma}[section]

\newtheorem{definition}{Definition}[section]

\newcommand{\mathsym}[1]{{}}

\def\pr{{\mathbf P}}

\title[Hidden Ancestor Graphs]{ 
Hidden Ancestor Graphs: Models for\\Detagging Property Graphs
}
\author[Darling, Clark]
{R.W.R.\ Darling \hspace{.6cm} Gregory S. Clark}
\address{National Security Agency, Fort George G Meade MD, USA}
\author[Tucker]{J. D. Tucker}
\address[SNL]{Sandia National Laboratories, Albuquerque NM, USA}
\date{\today}
\begin{document}
\maketitle

\begin{abstract}
Consider a graph $G$ where each vertex is visibly labelled as a member of a distinct class,
but also has a hidden binary state: \textit{wild} or \textit{tame}.
Edges with end points in the same class are called \textit{agreement edges}.
Premise: an edge connecting vertices in different classes -- a \textit{conflict edge} --
is allowed only when at least one end point is wild.
Interpret wild status as readiness to form connections with any other vertex,
regardless of class -- a form of class disaffiliation.
The learning goal is to classify
each vertex as wild or tame using its neighborhood data. In
applications such as communications metadata, bio-informatics, retailing, 
 or bibliography, adjacency in $G$ is typically created by paths of length two in a
transactional bipartite graph $B$. Class labelling, imported from a reference data source, is typically assortative, so agreement edges predominate. 
Conflict edges represent observed behavior (from $B$) inconsistent with
prior labelling of $V(G)$. Wild vertices are those whose label is uninformative.
The hidden ancestor graph constitutes a natural model for generating
agreement edges and conflict edges,
depending on a latent tree structure. The model is able to manifest
high clustering rates and heavy-tailed degree distributions
typical of social and spatial networks.
It can be fitted to graph data using a few measurable graph parameters, and
supplies a natural statistical classifier for wild versus tame.
\end{abstract}
{\small
\textbf{MSC class: } Primary: 05C80; Secondary: 90C35
}

\vspace{.4in}
\begin{center}
    \emph{We were born to be wild.}
\end{center}
\vspace{.15in}
\hspace{2.5in} --- Mars Bonfire, 1968
\vspace{.2in}

\newpage
\setcounter{tocdepth}{1}
\tableofcontents

\section{Introduction}
\subsection{Scope}
This paper models communities in graphs, but it is not about community detection.
It also deals with assortative vertex labels in graphs, but it is not
about label propagation. The context of this work is one familiar to graph data engineers,
 but largely ignored in the research literature. 
 It is one in which vertices of a large graph $G:=(V(G), E(G))$ arise as left nodes
 of a large bipartite graph $B$, expressing transactions between these left nodes and
 an even larger set of right nodes.
\begin{enumerate}
\item Adjacency in $E(G)$ is created by paths of length two in $B$, or maybe a subset of such paths subject to temporal constraints.
\item
     Labels are assigned to vertices in $V(G)$ using data from a source ouside $B$.
     \end{enumerate}
Usually end points of an edge $e \in E(G)$ have the same label, and we
call $e$ an \textit{agreement edge}. 
In the anomalous case  where end points have different labels we call $e$ a \textit{conflict edge}.
Conflict edges represent a discrepancy between observed behavior (from $B$)
and prior assumptions (labelling of $V(G)$).

\textbf{Premise: }\textit{
Each vertex has a hidden binary state, termed {\em tame} or {\em wild}.
Every conflict edge \textbf{must} have at least one wild end point.}

The premise is equivalent to the assertion that the wild vertices form a 
vertex cover of the conflict edges.
The interpretation of {\em tame} or {\em wild} is highly context-dependent:
see Table \ref{tab:divers} for examples. The term {\em wild} is not pejorative (think ``wild card''),
and indeed discovering the set of wild vertices is the operational goal.

\subsection{Concrete examples where hidden ancestor graph models may apply}\label{s:examples}
Table \ref{tab:divers} indicates four areas of application of our models.
 
 \begin{table}
    \centering   
    \begin{tabular}{lcccc}\toprule
    Field   &  Vertex & Label & Edge & Wild Vertex  \\ \midrule
    \textcolor{gray}{Retailing} &  retail & product & 2 items in same & item not specific  \\
  &  item & category & shopping basket & to one category \\ \midrule
 \textcolor{gray}{Telecoms} & network & service & 2 network IDs with & ID not indicative  \\
 & identifier & area  &$\geq 1$ common users & of user location   \\ \midrule
\textcolor{gray}{Informatics}  &  journal & subject & 
2 articles with & inter-disciplinary \\ 
      &  article & code & $\geq 1$ common authors &  article \\  \midrule
        \textcolor{gray}{Genomics} &  gene & phenotype & 2 mutations & not indicative   \\
   &  mutation  & &in same organism  & of phenotype  \\      
      \bottomrule
    \end{tabular}
\caption{Networks suitable for hidden ancestor graph models, and interpretation
of {\em wild vertex} in each. More details in Section \ref{s:examples}.
    }
    \label{tab:divers}
\end{table}

\subsubsection{Merchandise categorization for online retailing}
Suppose $V$ is a list of items offered for sale in an online marketplace.
Assume that a large sample $S$ of online transactions have been collected in the form of shopping baskets.
Connect two items by an edge of weight $w$ if they co-occur in
$w$ shopping baskets.
A conflict edge means a pair of items, in different product categories, appearing in the same shopping basket(s).
A wild vertex means an item which defies categorization, 
because it often appears in the same basket as seemingly unrelated items.

\subsubsection{Distributed communication networks}
Suppose $V$ is a list of points at which users gain access to a large communications network.
 Label these access points by intended service areas.
Business records for a limited time period 
are available as a temporal bipartite graph, where each link is a transaction between a user
(right node), and an access point (left node).
Insert an edge between two access points if at least one user touches both points within a short time period.
Agreement edges connect access points in the same service area. 
A conflict edge is a pair of access points, labelled with different service areas,
but which have common users.
A wild vertex means an access point not indicative of user location.

\subsubsection{Bibliographic data bases}
Suppose $V$ is a list of articles published in computer science journals in 2022.
Label each article by a primary code from the ACM Computing Classification System.
An edge $\{v, v'\}$ of weight $w$ between two articles $v$ and $v'$ is created
if there are $w$ authors in common between $v$ and $v'$.
A conflict edge is pair of articles with different codes, but at least one common author.
A wild vertex is a cross-disciplinary article. For a case study of this type, see
Section \ref{s:casestudy}.

\subsubsection{Genotypes with phenotype labels}
Suppose $V$ is a list of (sets\footnote{For brevity, we shall call a vertex a mutation, even when it is really a set of mutations.} of) gene mutations observed in organisms of some species, at some epoch.
Geneticists have attempted to label each mutation by a phenotype, meaning an observable physical
characteristic. Assume that a large sample $S$
of organisms is available. An edge $\{v, v'\}$ of weight $w$ between two mutations is created
if there are $w$ organisms in $S$ which exhibit both the mutations $v$ and $v'$.
Thus a conflict edge between two mutations occurs if they are associated with different phenotypes,
and yet occur in the same individuals. A wild vertex means a mutation with no specific phenotype expression.

\subsection{Statistical features of graphs we wish to model}
The class of graph models we wish to introduce often have the following properties:
\begin{enumerate}
\item[a.]
\textbf{Sparsity: }most vertices have low degree (i.e., number of connections)
\item[b.]
\textbf{Heavy tailed }vertex degree distribution, such as the Log-normal: a few vertices have very high degree. 
\item[c.]
High average \textbf{clustering }coefficient: if $v' \sim v$ and $v'' \sim v$, then
there is a good chance that $v' \sim v''$.
\item[d.]
\textbf{Assortative }labels: if $v \sim v'$, there is a good chance that $v$ and $v'$ have the same label.
\end{enumerate}
We encourage readers to convince themselves that properties c. and d. are likely to occur in all four 
examples listed in Section \ref{s:examples}. Properties a. and b. will also occur if large data sets are chosen.

There are two other properties we would like our models to have:
\begin{enumerate}
\item[e.]
\textbf{Parsimony: } A dozen parameters should suffice to control graph characteristics mentioned above.
\item[f.]
\textbf{Scalability: } Graphs from 1\% scale to 100\% scale should be fast and easy to construct, given sufficient memory.
\end{enumerate}

\subsection{Literature review}
\subsubsection{Hidden variable models}
The vast literature on complex networks contains a variety of generative models, 
which offer some but not all of the properties a. to f. mentioned above. 
Barth\'elemy \cite{bar} and \cite{barb} discuss {\em hidden variable models} for spatial networks,
where existence of an edge between two nodes is dependent on hidden variables at each node,
and possibly on distance between nodes if they are already embedded in a metric space.
We refer the reader to  \cite{bar} for references to half a dozen specific models, which show
that high clustering may occur despite overall sparsity of the graph, because
of existence of dense communities.
Our hidden ancestor graph model (Section \ref{s1:hag}) belongs to this class of hidden variable models, where
the hidden state of a vertex is its ancestry in a latent tree, and its wild or tame status.

\subsubsection{Models with a hidden vertex state}
Darling and Velednitsky \cite{vel} studied algorithms to detect a hidden wild state in
vertices of bipartite graphs. Both \cite{vel} and the present work fall under the heading of Section 2.2
(anomalies in static attributed graphs) of Akoglu, Tong, and Koutra's survey \cite{ako}
of graph anomaly detection. Section 2.2.2 of \cite{ako} is concerned with detecting outliers in
communities. In the methods summarized there, outlier detection is seen as part of community detection.
Such methods are not optimal for hidden ancestor graphs, where community affiliation is explicitly known.

\subsubsection{Variants of the stochastic block model}
The \textbf{stochastic block model}(SBM) \cite{sni} initially partitions vertices into color classes, and then prescribes a connection rate between each pair of vertices that depends only on the color class to which each belongs. 
Sparse models with a moderate number of color classes do not permit high clustering rates. 
In SBMs the statistician's task is to detect the communities from the graph
structure, and there is no hidden wild or tame vertex state. The
relationship between hidden ancestor graphs and SBMs is discussed in greater
depth in Section \ref{s:compareblockmod}. The analysis by Bickel and Chen \cite{bic},
repurposed in \cite{kar} for the degree corrected case, offers a template for deeper analysis of our model.
Generalizations of SBMs to multiplex networks are described by Barbillon et al \cite{barbie}.

\subsection{Contributions}

The contributions of the present work are three-fold:
\begin{enumerate}
    \item[(a)] 
 Define parametrized stochastic models\footnote{
 Hidden ancestor graphs were proposed in an earlier version \cite{hag}
of this work. Several years of experience led to the present
simplification of \cite{hag} and a robust fitting procedure.} 
applicable to graphs in Table \ref{tab:divers}, 
 and supply an explicit linear time construction of their edge set (Sections \ref{s1:bb}, \ref{s1:hag}).
     \item[(b)] 
Demonstrate the identification and fitting of model parameters (Sections
\ref{s1;hagparamid}, \ref{s1:haggen}).
    \item[(c)] 
Construct a vertex cover for the conflict edges (not necessarily a minimum vertex cover) 
which contains nearly all wild vertices and excludes nearly all tame vertices in such graph models
(Section \ref{s1:class}).
\end{enumerate}

\section{Basic building block of hidden ancestor graphs} \label{s1:bb}
\subsection{Triply randomized configuration model}\label{s:edgegenerator}
The configuration model is a standard stochastic technique for
generating a random graph having a prescribed degree
sequence \cite{fri}. Script \ref{s:edgegenscript} uses 
three levels of randomization to create multigraphs, eaach a variant of the Chung-Lu \cite{chu} random graph generation scheme.
The \textit{Mathematica} \cite{wol} code\footnote{
We ask the reader's forbearance for placing our scripts up front rather than
moving them to an Appendix or a web repository; 
the justification is that statistical dependencies and
computational efficiencies can be made precise in scripts, but not so easily in text.
Parse scripts using Language Reference \cite{wol}. Command line
\texttt{wolframscript} runs against the \textit{Mathematica} Engine without purchase of a Wolfram license.
} is especially concise.

\subsubsection{Mathematica script for \texttt{edgeGenerator}}\label{s:edgegenscript}

\begin{Verbatim}[frame=single,fontsize=\small]
edgeGenerator[v_, y_, t_] := Module[{boundaries, arrivals, vertexDegrees,
      vtx2Deg, halfedges}, 
  boundaries = Prepend[Accumulate[y], 0.0]; 
  arrivals = Drop[Drop[Flatten[RandomFunction[ PoissonProcess[t], 
      {0.0, Last[boundaries]}][[2,2]]], -1], 1]; 
  vertexDegrees = BinCounts[arrivals, {boundaries}]; 
  vtx2Deg = Select[Transpose[{v, vertexDegrees}], #1[[2]] > 0 & ]; 
  halfedges = Flatten[Table[Table[vtx2Deg[[j,1]], {vtx2Deg[[j,2]]}],
      {j, 1, Length[vtx2Deg]}]]; 
  Return[(UndirectedEdge[#1[[1]], #1[[2]]] & )
      /@ Sort /@ Partition[RandomSample[halfedges], 2]]; 
];     
\end{Verbatim}

\subsubsection{Parsing script \ref{s:edgegenscript}}
Our \texttt{edgeGenerator} function has arguments $\mathbf{v, y}, t$,
where $\mathbf{v}$ is a symbolic array of $n$ vertices, 
$\mathbf{y}$ is an array of $n$ positive numbers (such as
i.i.d. random variates from a Log-normal distribution), and $t$
is a positive number. We call $y_i$ the \textbf{mark} of vertex $v_i$,
and we call $t$ the \textbf{rate}. The function returns a list of
edges on vertices $\mathbf{v}$, which may include loops, and multiple
copies of the same edge. 

\textbf{Steps of edge generation}
\begin{enumerate}
    \item[(a)] 
    In \texttt{boundaries} and \texttt{arrivals}, split
    the interval from $0$ to $\sum y_i$ into $n$ intervals,
    of respective lengths $y_1, y_2, \ldots, y_n$. 
    Generate arrivals of a Poisson process of rate $t$ on
    the whole interval.
    \item[(b)]  
    In \texttt{vertexDegrees} and \texttt{vtx2Deg}, 
    associate to vertex $v_i$ the
    number $B_i$ of arrivals of the Poisson process in the $i$-th
    interval (whose length is $y_i$) omitting cases when there are no arrivals.    
    \item[(c)]  
    In \texttt{halfedges}, build a list which includes $B_i$ copies
    of vertex $v_i$, known as \textit{half-edges} in the configuration model.   
    \item[(d)]  
    In the \texttt{Return}, shuffle the list of half-edges,
    split them up into consecutive pairs, and convert these pairs
    into edges. It the length of the list is odd. the last half-edge
    is lost.
\end{enumerate}

\subsubsection{Interpreting the output}
The output of the algorithm may include loops, and multiple copies of the same edge.
In hidden ancestor graphs we discard the loops\footnote{In other contexts they could be used to model
self-interactions} to give a multigraph on vertices $\mathbf{v}$,
where some vertices may be isolated. When we use the term {\em valency}, we mean the
number of incident edges, excluding loops, not all of which need be distinct.
The sum of the valencies is denoted $2 |M|$, where $M$ is the edge list.
We shall use the term {\em neighbor count} to refer to the number of $j\neq i$
for which $v_j$ is connected to $v_i$ by at least one edge. 
The sum of the neighbor counts is denoted $2|E|$, where $E$ lists distinct edges.
Figure \ref{f:3graphssamemarks} illustrates these distinctions.

\subsubsection{Levels of randomization}
\begin{figure}[ht]
    \caption{\textit{Three instances (in red) of edge generation on twenty vertices
    with the same set of Log-normal marks. 
    Loops are omitted, but multi-edges are included,
    Edges of the two other graphs are shown in gray, so every gray edge is red in another graph.
    Thus the two vertices at extreme right are both present
    in the first graph, and both absent from the third graph.
    Five of the vertices were unlucky in that no incident
    edges appeared in any of the three graphs.}}
    \label{f:3graphssamemarks}
    	\begin{center}
\scalebox{0.42}
{\includegraphics{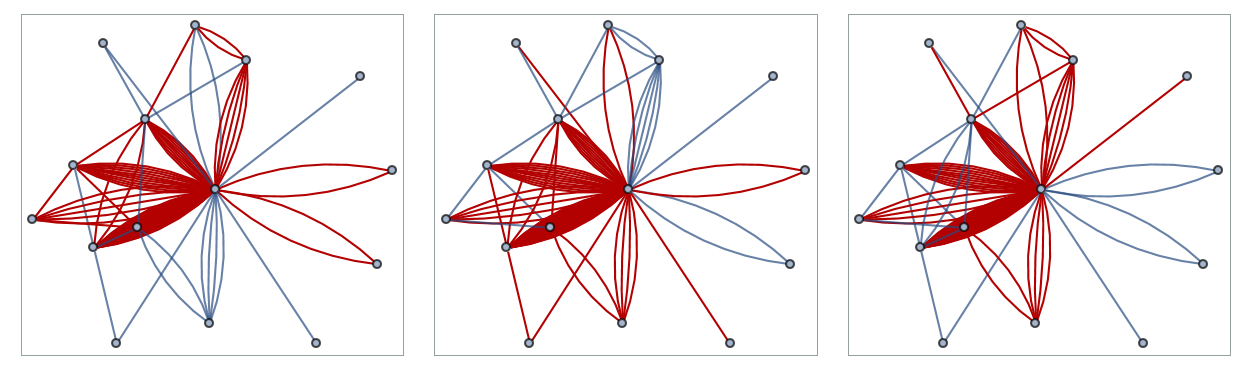}}
	\end{center}	
\end{figure}

The three levels of randomization in this generator are:
\begin{enumerate}
\item
Random generation of the marks $\mathbf{y}$.
\item
Generation of the Poisson process, which is a random
discretization of the marks.
\item
Shuffling of the half-edges to make the edges.
\end{enumerate}
Only the last of these occurs in the classical configuration model.
Figure \ref{f:3graphssamemarks} illustrates three random
graphs generated using the same 20 marks, to illustrate the
variation resulting from levels (2) and (3) above.
When isolated vertices are omitted from these plots,
it is important in the sequel to remember their existence,
because the denominator in statistics such as average vertex degree
is not the number of occupied vertices, but the number of vertices (and marks) supplied as input arguments.
We return to this point in Appendix \ref{a:censor}.

\subsection{Poisson Log-normal model}
Suppose the marks $(y_i)$ arise as i.i.d. Log-normal$(\phi, \sigma)$ random variates.
The half-edge counts follow the Poisson Log-normal model (Appendix \ref{a:quartilefit}),
used by Bulmer \cite{bul} to model species abundance. 
Log-normal leaf marks allow our edge generator to create graphs with heavy-tailed
vertex degree distribution (Figure \ref{f:valencies}). When the parameters $\phi$ and $\sigma$
are comparable in size, marks $y_i$ with value less than $\log{(2)}$ are common, with a likely 
result that vertex $v_i$ has no half-edges attached to it (Appendix \ref{a:censor}). 
Such a $v_i$ will be isolated.

\subsection{Estimating loop rates}
The number $B_i$ of half-edges for vertex $v_i$
is Poisson$(y_i)$, and $(B_1, B_2, \ldots, B_n)$ are mutually independent.
The total number $H:=\sum_i B_i$ of half edges is Poisson$(s)$, where $s:=\sum_i y_i$.
The conditional probability, given $(B_i)$, that the first edge produced in the construction above is a loop
on vertex $v_i$ is
\[
\frac{B_i (B_i - 1)}{H (H-1)}.
\]
Summing over $i$, the conditional probability that the first edge is a loop on some vertex is
\[
\frac{\sum_i B_i (B_i - 1)}{H (H-1)}.
\]
Elementary Poisson calculations show that the numerator of the last expression has mean $\sum_i y_i^2$,
and the denominator has mean $s^2$. If we adopt the heuristic of replacing the mean of a quotient of random
variables by the quotient of their means, we obtain a useful estimate:

\subsubsection{Heuristic}
\textit{For large $n$, the proportion of vertex pairs generated which are not loops is approximately}
\begin{equation}\label{e:non-collision}
    1 - \frac{\sum_i^n y_i^2}{(\sum_i^n y_i)^2}.
\end{equation}

In the special case where the $(y_i)$ arise as i.i.d. Log-normal$(\phi, \sigma)$ random variates,
the distribution of (\ref{e:non-collision}) depends only on $\sigma$, because factors of $e^{2 \phi}$
cancel in numerator and denominator. This fact will be exploited later.

\section{Hidden ancestor graph construction}\label{s1:hag}

\subsection{Vertex array: root, tribes, and clans}
\begin{figure}
    \caption{\textit{Latent tree with parameters $(b_0, b_1, b_2)=(4, 2, 3)$.
    The root in the center has
    four tribes, with two clans in each tribe, each of three siblings.
    The 24 vertices in the hidden ancestor graph are denoted $v[i, j, k]$.
    Tribes $t[i]$ and clans $c[i,j]$ are latent, not visible.
    }}
    \label{f:latenttree}
    	\begin{center}
\scalebox{0.35}
{\includegraphics{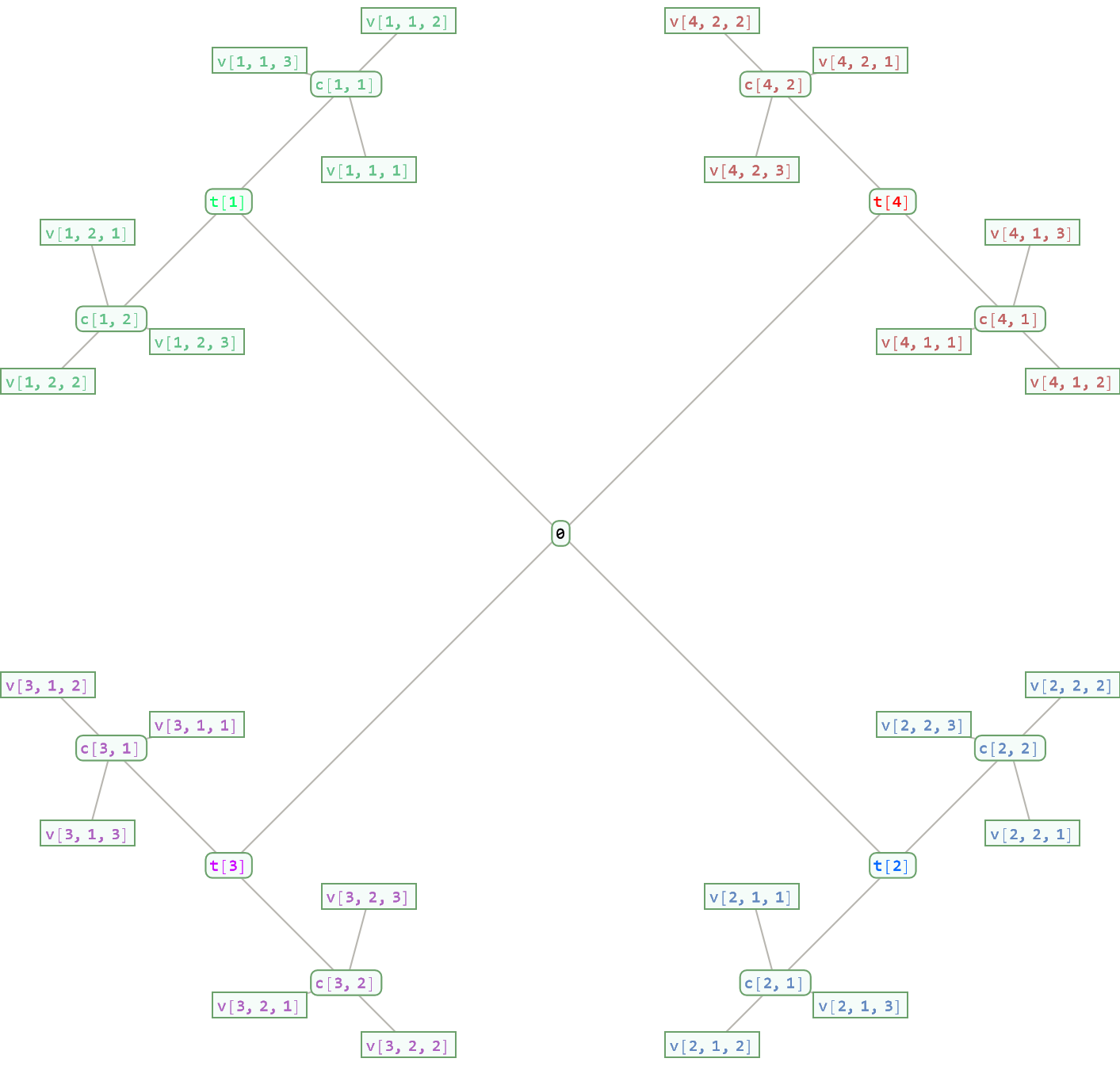}}
	\end{center}	
\end{figure}

We present the simplest non-trivial type of hidden ancestor graph, whose vertices
are the nodes at depth three\footnote{Depths four or more will be discussed in Section \ref{s:hagdepth4}.}
in a deterministic rooted tree.
The non-leaf nodes of this tree, and the {\em parent of} relations in the tree, do not form a visible part of
hidden ancestor graph, so we call it the \textbf{latent tree}. Its root is the 
eponymous \textbf{hidden ancestor}. 
The latent tree structure will be used to control the way in which edges of the 
hidden ancestor graph are generated.
Following branching process tradition we use the anthropogenic metaphor of \textit{tribes},
\textit{clans}, and \textit{siblings} as though the vertices of the graph were people, although Table \ref{tab:divers}
show that this is not an intended  application.
Use Figure \ref{f:latenttree} to visualize the following definition.

\begin{definition}[VERTEX ARRAY]\label{d:hagvert}
The vertex set of a (depth three) hidden ancestor graph is a triply indexed array 
\begin{equation}\label{e:vertexarray}
   \mathbf{v}:=\{v[i,j,k], 0 \leq i < b_0, 0 \leq j < b_1, 0 \leq k < b_2\}. 
\end{equation}
which refer to leaf vertices of a latent rooted tree: the ancestor has $b_0$ children, called \textbf{tribes};
each tribe has $b_1$ children, called \textbf{clans}; each clan consists of $b_2$ \textbf{siblings}.
Vertex $v[i,j,k]$ is $k$-th sibling in clan $c[i,j]$, which belongs to tribe $t[i]$.
\end{definition}

\begin{definition}[LABELS] \label{d:vertlabels}
Vertex $v[i,j,k]$ is labelled with the index $i$ of its tribe.
\end{definition}

Figure \ref{f:latenttree} shows vertex labels as colors. These colors or labels are observed
at least on those vertices with at least one incident edge. In other words, the tribe
to which a leaf of the latent tree belongs is revealed if the vertex is occupied.

\textbf{Example: }Vertices $v[1,3,2]$ and $v[1,2,2]$ are in the same tribe, but different clans. They carry the same label.

\subsubsection{Roles of parameters $(b_0, b_1, b_2)$}
Definition \ref{d:vertlabels} shows that $b_0$ is the number of vertex labels. Usually the clan size $b_2$ is less than a hundred, 
but even when it is
in the thousands simulation experiments show that average local clustering can take large values.
Use $b_1$ to control the size of the graph.
Statistical consistency is not studied here, but if it were
we would let $b_1 \to \infty$, holding $b_0$ and $b_2$ constant.

\subsubsection{Why are latent trees natural?}
Look at the examples of Table \ref{tab:divers}.\\
\textit{Genomics: }
Genetic mutations have a latent hierarchy based on genetic history of organisms.\\
\textit{Retailing: }
Retail hierarchy is based on industries, manufacturers, and suppliers.\\
\textit{Telecoms: }
Access points are leaf nodes of a control hierarchy of servers and routers \cite{sab}.\\
\textit{Bibliography: }
Articles have a latent hierarchical structure, in that several 2023 articles seek to answer a question posed in some article of a prior year.

\subsection{Hidden ancestor graph generation}
\subsubsection{Inputs}\label{s:haginputs}
Inputs to the hidden ancestor graph generator consist of 
\begin{enumerate}
    \item 
    the $b_0 \times b_1 \times b_2 $ vertex array described in Definition \ref{d:hagvert};
    \item 
     a $b_0 \times b_1 \times b_2 $ array of positive real numbers as vertex marks, called {\em leaf marks}, which are typically 
      i.i.d Log-normal$(\phi, \sigma)$ random variables;
        \item 
         a $b_0 \times b_1 \times b_2 $ boolean array which assigns wild (1) or tame (0) status (which will affect
         edge generation in Section \ref{s:hagedgegen}) to vertices,
    typically set by i.i.d Bernoulli$(\omega)$ random variables, for some $\omega \in (0,1)$;
    \item 
    a probability vector $(q_0, q_1, q_2)$.
\end{enumerate}
Examples of parameter choices used our illustrations appear in Table \ref{tab:genericparameters}.
While graphical illustrations necessarily show small graphs, in
industrial use cases hidden ancestor graphs may have millions
of vertices and hundreds of millions of edges.

\begin{table}[]
    \centering
    \begin{tabular}{c|c|c|c|c|c|c|c|c}
    \toprule
         $b_0$ & $b_1$ & $b_2$ & $q_0$ & $q_1$ & $q_2$ & $\phi$ & $\sigma$ &$\omega $\\  \midrule
    80 & 50 & 25 & 0.19 & 0.1 & 0.71 & 1.79 & 2.15 & 0.042 \\
\bottomrule
    \end{tabular}
    \caption{Typical input parameters for
    Section \ref{s:haginputs}. Leaf marks are i.i.d. Log-normal$(\phi, \sigma)$;
    wild status is i.i.d. Bernoulli$(\omega)$. Graph instances have about 1.6M edges 
    (after removing 38\% which formed loops) on 100K vertices,
    of which about 81K are visible in the agreement graph. The conflict graph has a giant component.
    Observables $(d_C, p_C, d_A, \kappa_A)$ typically $(0.7, 0.16, 30, 0.43)$.    
    }
    \label{tab:genericparameters}
\end{table}

\subsubsection{Outputs}
The script \ref{sc:hag} produces three (possibly overlapping) lists of edges,
which may include loops. Subsequently these lists are concatenated, and loops are removed,
creating a multigraph. Although the script presents vertices in $v[i,j,k]$ notation, the
user should suppose that vertex names have been hashed, obscuring the ancestry in the latent tree,
which survives only in the form of the (tribal) labelling function
\begin{equation}
    v[i,j,k] \mapsto i \in \{0, 1, \ldots,b_0-1\}.
\end{equation}
The following are \textbf{not} provided in the output:
\begin{enumerate}
    \item[(a)] The leaf marks.
    \item[(b)] The wild or tame status of vertices.
    \item[(c)] Latent ancestry of vertices (except for tribal label).
    \item[(d)] Vertices in the list (\ref{e:vertexarray}) not incident to any edge.
\end{enumerate}
The number of isolated vertices can be known by subtraction if the product $n:=b_0 b_1 b_2$ is known, 
or possibly estimated from quartiles as in Appendix \ref{a:quartilefit}.

\subsubsection{Two edge types}
\begin{figure}
    \caption{\textit{Agreement and conflict graphs taken from a single hidden ancestor graph instance
    with $b_0 = 6, b_1 = 4, b_2 = 20$. Each edge on the left has end points of the same color (tribe);
    on the right, end points of different colors. Neither graph need be connected.
    Out of 480 vertices, 387 have at least one agreement edge, 109 at least one
    conflict edge, and 93 were isolated. Most vertices in the conflict graph
    are tame vertices with a wild neighbor.}
    }
    \label{f:agreeplusconflict}
    	\begin{center}
\scalebox{0.45}
{\includegraphics{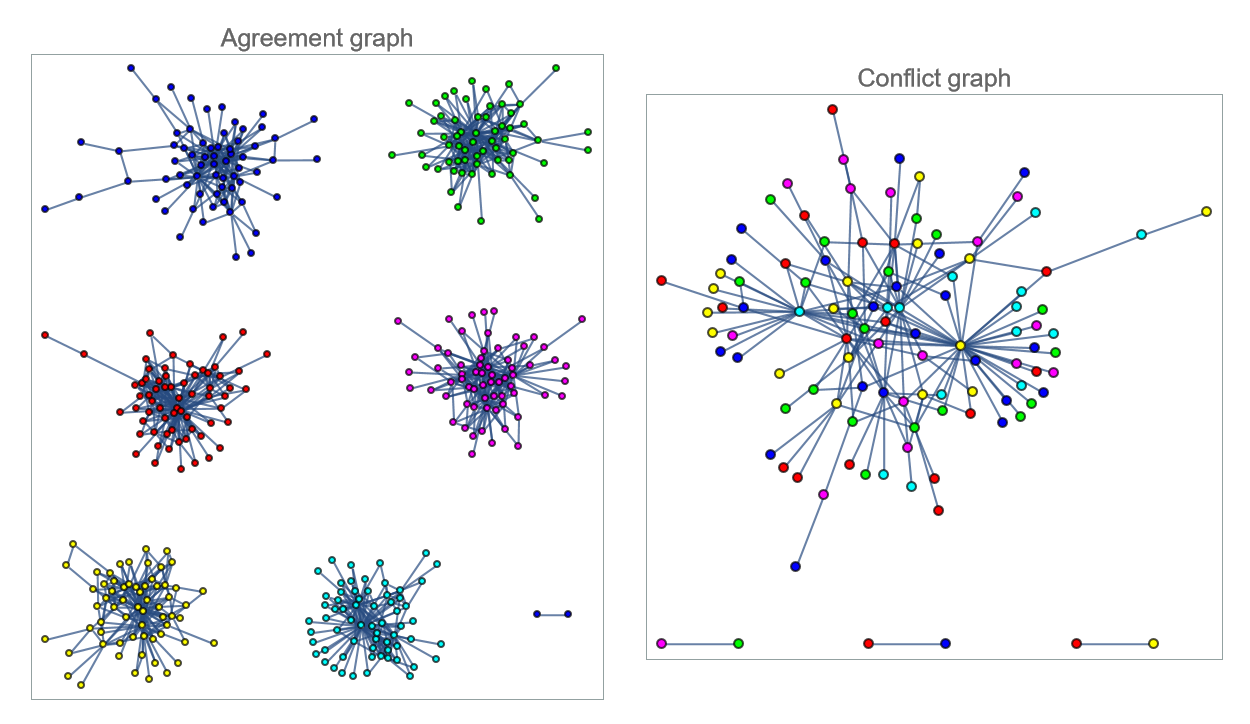}}
	\end{center}	
\end{figure}

\begin{definition}[EDGE TYPES]\label{d:edgetypes}
The hidden ancestor (multi)graph admits two types of edges, depending on labels of the end points of the edge:
\begin{enumerate}
    \item[agreement: ]end points have the same label, i.e. leaves in the same tribe of the latent tree.
    \item [conflict: ]different labels (end points leaves are in different tribes)
    -- but this is allowed \textbf{only} if at least one end point has wild status\footnote{
    No edge is allowed between leaf vertices in different tribes of the latent tree if both vertices are tame.}.
\end{enumerate}
The induced subgraphs (which usually have vertices in common) are called the agreement graph and
conflict graph, respectively. A vertex with at least one incident conflict edge is called a conflict vertex.
\end{definition}
Figure \ref{f:agreeplusconflict} illustrates the agreement graph and conflict graph;
conflict vertices typically have incident agreement edges.

\subsubsection{Script for \texttt{hag3Edges}} \label{sc:hag}

\begin{Verbatim}[frame=single,fontsize=\small]
hag3Edges[vertices_, marks_, wildStatus_, q0_, q1_, q2_] := Module[
        {b0, b1, b2, wildMap, preedges0, edges0, edges1, edges2}, 
   {b0, b1, b2} = Dimensions[marks]; 
   wildMap =  Association[(#1[[1]] -> #1[[2]] &) /@
       Transpose[{Flatten[vertices], Flatten[wildStatus]}]]; 
   preedges0 = edgeGenerator[Flatten[vertices], Flatten[marks], q0]; 
   edges0 = Select[preedges0, #1[[1, 1]] == #1[[2, 1]] ||
        Lookup[wildMap, #1[[1]]] + Lookup[wildMap, #1[[2]]] > 0 &]; 
   edges1 = Flatten[Table[edgeGenerator[
       Flatten[vertices[[j]]], Flatten[marks[[j]]], q1], {j, 1, b0}]]; 
   edges2 = Flatten[Table[edgeGenerator[vertices[[j, k]],
       marks[[j, k]], q2], {j, 1, b0}, {k, 1, b1}]]; 
   Return[{edges0, edges1, edges2}];
];    
\end{Verbatim}

\subsubsection{Parsing the script \ref{sc:hag}}\label{s:hagedgegen}
\begin{enumerate}
    \item[(0)] 
    First apply the \texttt{edgeGenerator} of Section \ref{s:edgegenscript}
    to the flattened array of all $n:=b_0 b_1 b_2$
    vertices, with all $n$ marks, with rate $q_0$. However retain an edge only if its endpoints belong to the same tribe,
    or else if at least one end point is wild: these appear as \texttt{edges0} of the script \ref{sc:hag}. Call these
    \textbf{root edges}; most are conflict edges.
    \item[(1)] 
    Second apply the \texttt{edgeGenerator} $b_0$ times, to vertices in each tribe, respectively,
    using the $b_1 b_2$ marks associated with that tribe, with rate $q_1$. Since both end points of such edges are in the same
    tribe they are all agreement edges, listed in \texttt{edges1} of the script \ref{sc:hag}.
    Call these \textbf{tribal edges}.
    \item[(2)] 
    Finally apply the \texttt{edgeGenerator} $b_0 b_1$ times, to vertices in each clan, respectively,
    using the $b_2$ marks associated with that clan, with rate $q_2$. These are also agreement edges, listed in \texttt{edges2} of the script \ref{sc:hag}. Call these \textbf{clan edges}.
\end{enumerate}
The script returns the root, tribal, and clan edges separately, with loops included.
There could be overlap among the three sets of edges.
Loops are stripped out later.

\begin{definition}[EDGE SET]\label{d:hagedges}
The edge set of the hidden ancestor (multi)graph is the concatenation of lists of root edges, tribal edges, 
and clan edges, excluding loops. Multiple instances of the same edge are retained.
\end{definition}

\subsection{Hidden ancestor graphs versus stochastic block models}\label{s:compareblockmod}
Suppose leaf marks take a constant deterministic value $t$, all vertices have wild status, $b_2 = 1$,
and $q_2 = 0$. This gives a form of stochastic block model \cite{bic, kar,sni} with $b_0$ blocks, and $b_1$ vertices in each block.
We expect a total of $q_0 b_0 b_1 t/2$ edges among all vertices, together with $q_1 b_1 t/2$ additional edges
among each of the $b_0$ blocks. In such models, tribal labels are not displayed, and the statistician is
expected to infer them from the graph community structure. In a large, sparse, stochastic block model,
average local clustering is insignificant.

The hidden ancestor graph goes beyond the stochastic block model in several ways:
\begin{enumerate}
    \item Hidden wild or tame status of vertices, with application to anomaly detection, not community detection.
    \item Possibility of large average local clustering when $b_2$ is comparable to mean agreement degree, even if the graph is sparse.
    \item Heavy-tailed edge multiplicities and vertex degrees, by suitable choice of the probability law of the marks.
\end{enumerate}

\subsection{Heavy-tailed distributions} 
\begin{figure}
    \caption{\textit{Valency refers to number of edges incident to a vertex, including 
    multiple copies of the same edge. The different shapes for conflict versus agreement
    edges reflect the censoring of loops among proposed clan edges, which make up
    the majority of agreement edges in a high clustering case.
    }
    }
    \label{f:valencies}
    	\begin{center}
\scalebox{0.47}
{\includegraphics{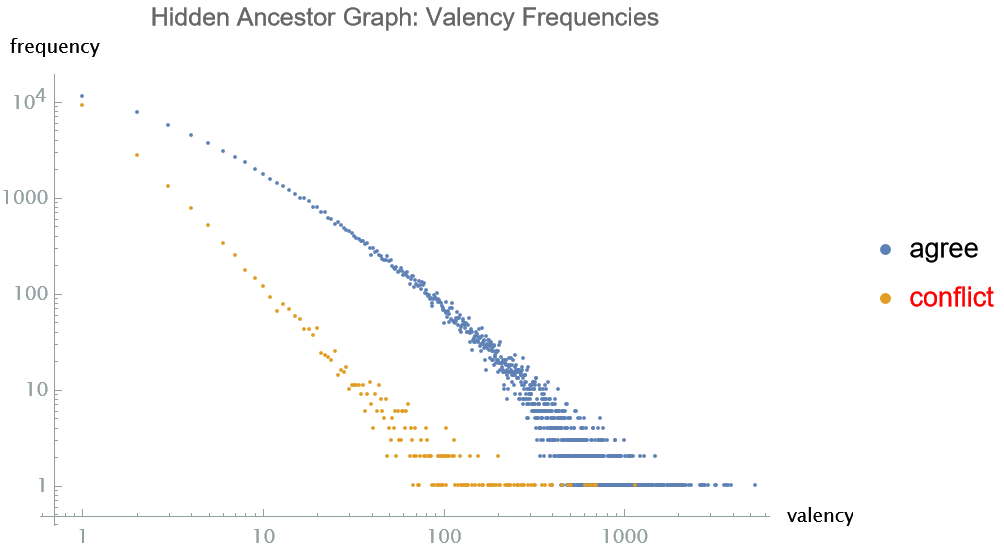}}
	\end{center}	
\end{figure}

See Mitzenmacher \cite{mit} for a discussion of why heavy-tailed distributions emerge in
large graphs, and Antoniou et al \cite{ant} and Radicchi et al \cite{rad} for instances
of the Log-normal.
Choice of the probability distribution of the marks array supplied to the script \ref{sc:hag}
will be reflected in vertex valencies, and multiplicities at each edge position.
Figures \ref{f:valencies} and Figure \ref{f:heavytail}
refer to a graph with parameters from Table \ref{tab:genericparameters}.

\subsubsection{Vertex valencies in the agreement and conflict graphs}

The nearly linear plot in log-log scale
illustrates the heavy-tailed valencies resulting from choice
of i.i.d. Log-normal leaf marks.

\subsubsection{Edge multiplicities}\label{s:multiplicities}
\begin{figure}
    \caption{\textit{Horizontal axis shows values of agreement edge multiplicity;
    for example 255K out of
    of 399K edge positions had multiplicity 1 (top left). 
    Vertical axis shows the number of positions for
    which that multiplicity occurs,
     for a hidden ancestor graph with Log-Normal distribution of marks.
    Some edge positions have multiplicity between 100 to 1000 (bottom right).
    The log-log scale displays a heavy-tailed pattern typical of
    large communication networks \cite{ant}.}
    }
    \label{f:heavytail}
    	\begin{center}
\scalebox{0.5
}
{\includegraphics{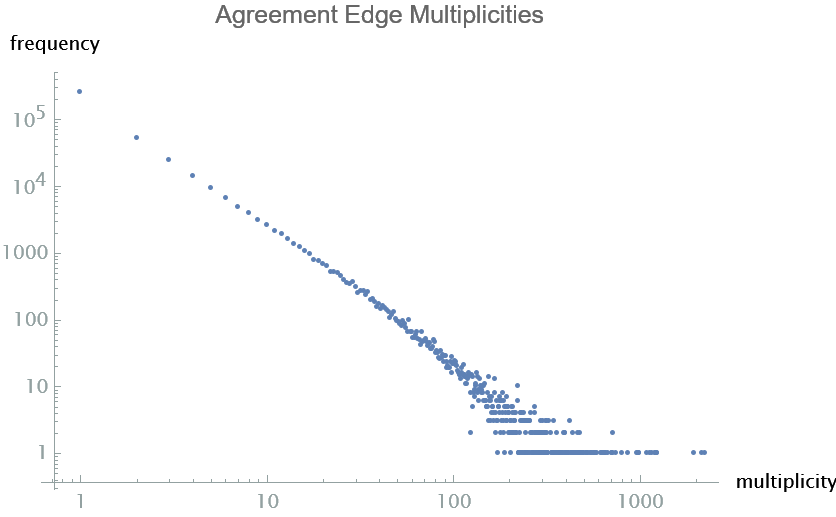}}
	\end{center}	
\end{figure}

 When a large mark occurs on a vertex, 71\% of its weight
is directed towards making connections with the 24 siblings of that vertex; so 
edge positions of large multiplicity are likely. The mean multiplicity of agreement edges
is 3.8 in the example of Figure \ref{f:heavytail}. By contrast,
multiplicity is nearly always 1 in the conflict graph.

\subsection{Observable statistics of the hidden ancestor graph}\label{s:observables}
Suppose that $b_0, b_1, b_2$ are known, hence so is the total vertex count $n:=b_0 b_1 b_2$.
The statistician does not see the latent tree, nor the hidden wild or tame status of vertices,
but at least can see vertex labels, and therefore can assemble conflict edges into the \textbf{conflict graph}
and agreement edges into the \textbf{agreement graph}. Here are the four principal observables:
\subsubsection{Mean agreement degree $d_A$:}
 twice the number of agreement edges, divided by $n$; the denominator
is usually larger than the union of vertices appearing in the edge set.
\subsubsection{Mean conflict degree $d_C$: }
twice the number of conflict edges, divided by $n$.
\subsubsection{Proportion $p_C$ of vertices in the conflict graph: }number of vertices with at
least one incident conflict edge, divided by $n$.
\subsubsection{Average local clustering coefficient $\kappa_A$ in the agreement graph: }
Local clustering at a vertex is the
proportion of the number of links between the vertices within its neighbourhood divided by the 
number of links that could possibly exist between them, taking no account of edge multiplicity.
Implicitly the denominator in the average is the number of vertices with at
least one incident agreement edge. We ignore clustering in
the conflict graph, which is usually negligible.

\subsection{Numerical example with output statistics}
Table \ref{tab:hag900K} shows an example of a hidden ancestor graph created by the {\em Mathematica} script \ref{sc:hag}
with $0.9$M. vertices (83\% of which were occupied), and $18.1$M. edges, of which $5.54$M. were distinct.
Coding languages such as Java\copyright \, can create examples fifty times larger
\cite{jhag}.
The first pair of rows shows input parameters. The second row shows the four observable output parameters above.

\begin{table}[]
    \centering
    \begin{tabular}{c|c|c|c|c|c|c|c|c}
    \toprule
        INPUT & $b_0$ & $b_1$ & $b_2$ & $q_0$ & $q_1$ & $q_2$ & $\phi$ & $\omega $\\ 
   PARAMETERS      & 200 & 150 & 30 & 0.15 &  0.13 & 0.72 & 1.95 & 0.04\\ \midrule
   GRAPH & $n$ & $|M|$ & $|E|$ & $\frac{|V|}{n}$ & $d_C$ & $p_C$ & $d_A$ & $\kappa_A$ \\
        STATISTICS & 0.9M & 18.1 M & 5.54 M & 0.83 & 0.8078 & 0.1655 & 39.39 & 0.4015 \\ \midrule
   ESTIMATES & \, & \, & \, & $\hat{q}_0$ & $\hat{q}_1$ & $\hat{q}_2$ & $\hat{\phi}$ & $\hat{\omega}$\\   
  from $(d_C, p_C, d_A, \kappa_A)$  & \, & \, & \, & 0.124 & 0.125 & 0.751 & 1.93 & 0.048 \\ \bottomrule
    \end{tabular}
    \caption{ For inputs in rows 1 and 2,
    see Section \ref{s:haginputs}; marks are i.i.d. Log-normal$(\phi, \sigma)$, taking $\sigma = 2.1$. 
    Rows 3 and 4 refer to a hidden ancestor graph generated with these inputs.
    $M$ refers to all edges (excl. loops); $E$ refers to
    filled edge positions; $V$ refers to occupied vertices; Section \ref{s:observables}
    defines observables $(d_C, p_C, d_A, \kappa_A)$.
   Rows 5 and 6 refer to estimates of parameters from observables 
using methods of Section \ref{s1;hagparamid} (step size $\delta=0.002$). Deviation between truth and estimate of 20\% is typical.
    }
    \label{tab:hag900K}
\end{table}

\section{Hidden ancestor graph parameters identification: simple case}\label{s1;hagparamid}

\subsection{Identifiability in statistics}
The notion of parameter identifiability in statistics is explained in Basu \cite{bas}.
The identifiability issue for hidden ancestor graphs arises from a hypothetical scenario where the
{\em Customer} visits the {\em Statistician} and asks her to model a large multigraph, each of whose vertices
carry one of $T$ visible labels. The number of agreement edges and conflict edges can be counted, as can
the number of conflict vertices. Vertex degree distribution appears to have a Log-normal tail.
Average local clustering coefficient in the agreement graph is also available,
possibly by sampling if the graph is very large. 

The {\em Statistician} wants to try to fit a hidden ancestor graph model with Log-normal leaf marks.
Is there a set of parameter choices which will serve this customer? If so, is there exactly one set;
in other words, are the parameters identifiable?

The problem is not yet well-posed. What exactly are we allowed to assume, and what remains to be estimated?
As a first step, we state a simplified version and solve it empirically. The full version is covered in 
Section \ref{s1:haggen}.

\subsection{Simplified parameter identification problem}\label{s:simpleidentify}
Consider a game whose players are called {\em Modeller} and {\em Statistician}. The 
{\em Modeller} creates a hidden ancestor graph with i.i.d. Log-normal$(\phi, \sigma)$ leaf marks
using a parameter set such as the one
listed in rows 1 and 2 of Table \ref{tab:hag900K}. Then he passes the graph edges to the 
{\em Statistician}, together with values of parameters $(b_0, b_1, b_2)$ and $\sigma$;
however the $\phi$ parameter is \textbf{not} shared.

The {\em Statistician} may now compute the four observables shown in rows 3 and 4 of Table \ref{tab:hag900K}
\[
(d_C, p_C, d_A, \kappa_A)
\]
since she knows the correct denominator $n:=b_0 b_1 b_2$ for the first three of these.
The game is: the {\em Statistician} must guess the other parameters
$(\omega, \phi, q_0, q_1, q_2)$ concealed by the {\em Modeller}. She can ask
the {\em Modeller} to generate as many instances of hidden ancestor graphs with the same parameters as she wants.

\begin{definition}\label{d:identify}
    The hidden ancestor graph parameter identification problem is: 
    determine whether there is only one
    probability vector $(q_0, q_1, q_2)$, wildness rate $\omega$, and $\phi$ parameter 
    which produces given $(d_C, p_C, d_A, \kappa_A)$, averaged over many instances of graph generation.
    If so, what is $(\omega, \phi, q_0, q_1, q_2)$?
\end{definition}
The rest of this section is devoted to describing a strategy for the {\em Statistician}.
This leads to the conclusion, which is not formally proved, that $(q_0, q_1, q_2)$, $\omega$, and $\phi$ are identifiable
under the given assumptions. We revisit this issue later in Section \ref{s:consistency}.

\subsection{Separation of variables}
The separation of the fitting problem into conflict and agreement parts is possible in
the Log-normal case because of the following Lemma, whose proof
follows immediately from the multiplicative role of $e^{\phi}$ in Log-normal$(\phi, \sigma)$ variates.

\begin{lemma}\label{l:tribalsplit}
Consider a hidden ancestor graph with i.i.d. Log-normal$(\phi, \sigma)$ leaf marks, 
and probability vector $(q_0, q_1, q_2)$. 
\begin{enumerate}
    \item[(a)] 
The root edges have the same law as those arising from Log-normal$(\phi_0, \sigma)$ leaf marks
and probability vector $(1, 0, 0)$, where $e^{\phi_0}:= q_0 e^{\phi}$.
    \item[(b)] 
Tribal and clan edges combined have the same law
those arising from Log-normal$(\phi', \sigma)$ leaf marks
and probability vector $(0, q_1', q_2')$, where 
\[
e^{\phi'} = (q_1 + q_2) e^{\phi}; \quad q_1' = \frac{q_1}{q_1 + q_2}  \quad q_2'=\frac{q_2}{q_1 + q_2}.
\]
\end{enumerate}
\end{lemma}

\begin{figure}
    \caption{\textit{Logical dependencies among 7 subroutines for parameter identification}}
    \label{f:dependencies}
    	\begin{center}
\scalebox{0.45}
{\includegraphics{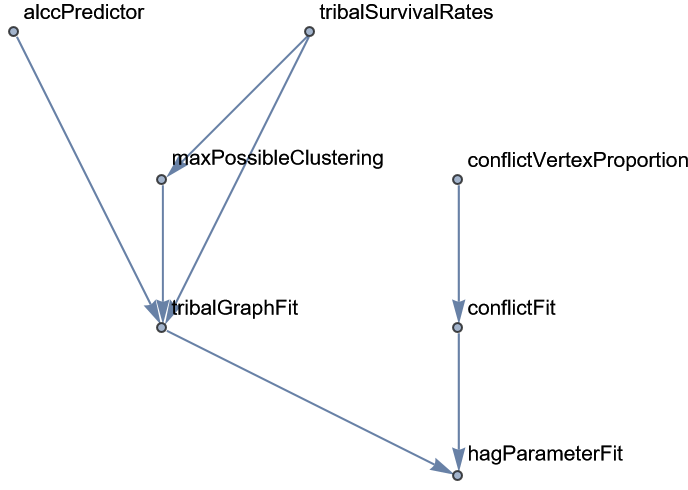}}
	\end{center}	
\end{figure}

Using a change of variables from $(\phi, q_0, q_1, q_2)$ to $(\phi_0, \phi_1, \phi_2)$ where
\[
(e^{\phi_0}, e^{\phi_1}, e^{\phi_2}):= (q_0 e^{\phi}, q_1 e^{\phi}, q_2 e^{\phi})
\]
we can separate the function
\[
(\omega, \phi, q_0, q_1, q_2) \mapsto (d_C, p_C, d_A, \kappa_A)
\]
into a pair of functions:
\begin{equation}\label{e:sepvar}
    (\omega, \phi_0) \overset{f_C} \longrightarrow (d_C, p_C) \quad
    (\phi_1, \phi_2) \overset{f_A} \longrightarrow (d_A, \kappa_A).
\end{equation}
Inversion of $f_C$ and $f_A$ are the topics of Sections \ref{s:conflgraphfit}
and \ref{s:agreegraphfit}, respectively.
Subroutines for inverting $f_C$ and $f_A$ are shown in Figure \ref{f:dependencies}.
Two nodes on the right refer to subroutines which together fit $\omega$ and $\phi + \log{q_0}$ from
the statistics $d_C, p_C$. The four nodes at upper left 
refer to subroutines which together fit $\phi + \log{q_1}$ and $\phi + \log{q_2}$ from $d_A, \kappa_A$.
The node at the bottom combines all this to fit $(\omega, \phi, q_0, q_1, q_2)$.
Figures in this section and the next are based on hidden ancestor graph
instances with parameters shown in Table \ref{tab:genericparameters}.

We will not indulge in the luxury of generating many graph instances, but will content
ourselves with fitting $(\omega, \phi, q_0, q_1, q_2)$ from a sample of size one. 
We present actual scripts for the fitting process,
rather than a mere narrative description;
our approach is more precise and takes up no more space.

\subsection{Fitting the conflict graph}\label{s:conflgraphfit}


\subsubsection{Mathematica script for \texttt{conflictVertexProportion}}\label{sc:cvp}

\begin{Verbatim}[frame=single,fontsize=\tiny]
conflictVertexProportion[b0_, b1_, b2_, Phi0_, Sigma_, Omega_] := Module[{y, z, tameP, cP}, 
    y = RandomVariate[LogNormalDistribution[Phi0, Sigma], b0*b1*b2]; 
    z = Select[(RandomVariate[PoissonDistribution[#1]] & ) /@ y, #1 > 0 & ]; 
    tameP = 1. - Omega ;
    cP = (1. - tameP^(#1 + 1) & ) /@ z; 
    Return[N[Total[cP]/(b0*b1*b2)]]; 
]; 
\end{Verbatim}
\subsubsection{Parsing the script \ref{sc:cvp} }
Here $\phi_0$ refers to $\phi + \log{q_0}$, so that Log-normal$(\phi_0, \sigma)$ 
is the distribution of the leaf mark portion
 devoted to root level edge generation as in Lemma \ref{l:tribalsplit}. 
 Generate these partial marks as the vector $\mathbf{y}$.
The corresponding half-edge counts (when non-zero) give the vector $\mathbf{z}$ in the script, indexed by vertices
incident to at least one root level edge. We claim the chance that such a vertex $v$, with $b \geq 1$ root level half-edg
is in the conflict graph is
\[
\omega + (1 - \omega) (1 - (1 - \omega)^b) = 1 - (1 - \omega)^{b+1}.
\]
This formula, expressed in the \texttt{cP} line, holds because $v$ is itself wild with probability $\omega$,
and if tame it is still in the conflict graph unless all its $b$ half-edges team up with tame vertices.
The output of the function is a Monte Carlo estimate of the proportion in the conflict graph.

\begin{figure}
    \caption{\textit{Fitting the conflict graph: the wildness rate $\omega$
    ($x$-axis) IS incremented until the implied proportion of
    vertices in the conflict graph ($y$-axis) attains the desired level,
    shown as a horizontal line. Noise is evident since 
    \texttt{conflictVertexProportion} samples Poisson Log-normal variates.
}}
    \label{f:wildfit}
    	\begin{center}
\scalebox{0.4}
{\includegraphics{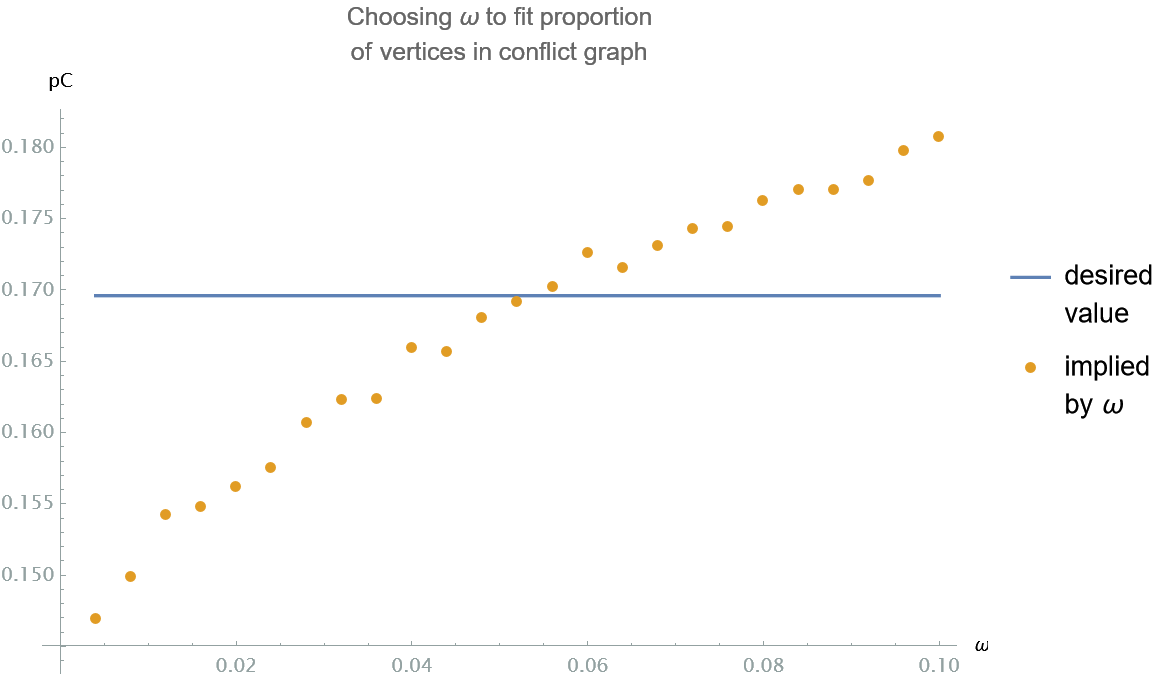}}
	\end{center}	
\end{figure}

\subsubsection{Mathematica script for \texttt{conflictFit}}\label{sc:confit}
\begin{Verbatim}[frame=single,fontsize=\small]
conflictFit[b0_, b1_, b2_, Sigma_, dC_, pC_, Delta_] := Module[{Omega, f, cvp}, 
    f[w_] := Log[dC*(b0/(w*(2. - w)*(b0 - 1)))] - Sigma*(Sigma/2); 
    {Omega, cvp} = {0., 0.}; 
    While[cvp < pC && Omega < 1., 
     Omega = Omega + Delta; 
     cvp = conflictVertexProportion[b0, b1, b2, f[Omega], Sigma, Omega]; 
    ]; 
   Return[{Omega, f[Omega]}];
]; 
\end{Verbatim}
\subsubsection{Parsing the script \ref{sc:confit}}
Inputs include the observables $d_C, p_C$ and a step size $\delta$.
The \texttt{While} loop increases $\omega$ by $\delta$ at each step, and uses Script \ref{sc:cvp}
to evaluate the resulting proportion in the conflict graph, stopping when $p_C$ is exceeded, as shown in Figure \ref{f:wildfit}.
The function \texttt{f} is designed to supply a value of $\phi_0$ so that
the expected number of conflict edges per vertex takes the intended value:
\[
\omega (2 - \omega)  \left(1 - \frac{1}{b_0}\right) \exp{(\phi_0 + \sigma^2)} = d_C.
\]
The mean of a Log-normal$(\phi_0, \sigma)$ is multiplied by the chance $1 - (1-\omega)^2$ 
that at least one end points of an edge is wild, and
by the chance the two end points are in different tribes. The function returns a pair $(\omega, \phi_0)$, which completes the fitting of the conflict graph. Interpret $\phi_0$ as the same as $q_0 e^{\phi}$ in the sense
of Lemma \ref{l:tribalsplit}.
\subsection{Three subroutines for agreement graph fitting}
This section describes the three subroutines at the top left of Figure \ref{f:dependencies}. 

\subsubsection{Mathematica script for \texttt{tribalSurvivalRates}} \label{sc:tribsurv}
The function \texttt{loopAvoidRate} with a vector argument encodes
formula (\ref{e:non-collision}).
\begin{Verbatim}[frame=single,fontsize=\small]
tribalSurvivalRates[b0_,b1_,b2_,Sigma_]:=Module[{y, yByTribe,means1,weights1, 
      survivalRates1 ,yByClan, means2,weights2, survivalRates2,Rho1, Rho2},
  y = RandomVariate[LogNormalDistribution[1.0, Sigma], {b0, b1,b2}]; 
  yByTribe = Table[Flatten[y[[j]] ],{j,1,b0}];
  yByClan=Flatten[Table[y[[j, k]],{j,1,b0},{k,1,b1}], 1];
  {survivalRates1, survivalRates2} = {Map[loopAvoidRate,yByTribe], 
     Map[loopAvoidRate, yByClan]};
  {means1,means2 } ={ Map[Mean,  yByTribe], Map[Mean, yByClan]};
  {weights1,weights2} = {means1/Total[means1],means2/Total[means2]};
  {Rho1, Rho2}= {Dot[weights1, survivalRates1 ],Dot[weights2,survivalRates2] };
  Return[{Rho1, Rho2}]; 
];
\end{Verbatim}
\subsubsection{Parsing the script \ref{sc:tribsurv} }
When a hidden ancestor graph is generated using leaf marks with a heavy-tailed
distribution, a third to a half of clan edges are commonly discarded because
of loops. To fit the agreement graph correctly,  start by
estimating the proportion of the mark which is converted to loop-free half-edges
at the tribal and clan level. These are  the outputs $(\rho_1, \rho_2)$
of the Script \ref{sc:tribsurv}. 
When leaf marks are Log-normal$(\phi, \sigma)$ random variates,the $\phi$
plays no part, and is set to the arbitrary value of 1.0  in the script.
The script generates a complete set of such leaf marks,
and uses heuristic (\ref{e:non-collision}) as the \texttt{loopAvoidRate}
function, and tribal and clan levels, respectively.

\subsubsection{Mathematica script for \texttt{maxPossibleClustering}}\label{sc:maxalcc}
The \texttt{loopRemover} function merely selects edges whose end points differ.
\begin{Verbatim}[frame=single,fontsize=\small]
maxPossibleClustering[b0_,b1_,b2_,Sigma_, dA_]:=Module[
           {Rho2,Phi2,tribalVertices ,tribalMarks,edges2},
  Rho2 = tribalSurvivalRates[b0,b1,b2,Sigma][[2]];
  Phi2 = Log[dA/Rho2] - Sigma*Sigma/2.0;
  tribalVertices = Array[v, {b1, b2}];
  tribalMarks = RandomVariate[LogNormalDistribution[Phi2, Sigma], {b1, b2}];
  edges2=loopRemover[Flatten[Table[
        edgeGenerator[tribalVertices[[k]], tribalMarks[[k]],1.0], {k,1,b1}]]];
  Return[N[MeanClusteringCoefficient[Graph[edges2]]]];
];
\end{Verbatim}

\subsubsection{Parsing the script \ref{sc:maxalcc} }
We seek the maximum value $\kappa_{\max}$ of the average local clustering coefficient (ALCC)
for a given value of $d_A$, if all the agreement edges in the tribal subgraph are clan edges.
If this maximum value were less than the observed value of $\kappa_A$, then it
would be impossible to fit a hidden ancestor graph model without reducing the value of $b_2$.
In the script \ref{sc:maxalcc}, we use the \texttt{tribalSurvivalRates} script to determine $\rho_2$,
and pick $\phi_2$ so that the Log-normal$(\phi_2, \sigma)$ distribution (for clan edges) has a mean
of $d_A/\rho_2$. Next generate a graph consisting of $b_1$ clan subgraphs, and return its ALCC $\kappa_{\max}$.

\begin{figure}
    \caption{\textit{Average local clustering coefficient (ALCC) in the agreement graph, as a
    function of loading on to clan (rather than tribal) edges.
    The horizontal line shows the desired ALCC. These Monte Carlo estimates 
    display high noise, making the intersection point imprecise.
    }}
    \label{f:clustering}
    	\begin{center}
\scalebox{0.375}
{\includegraphics{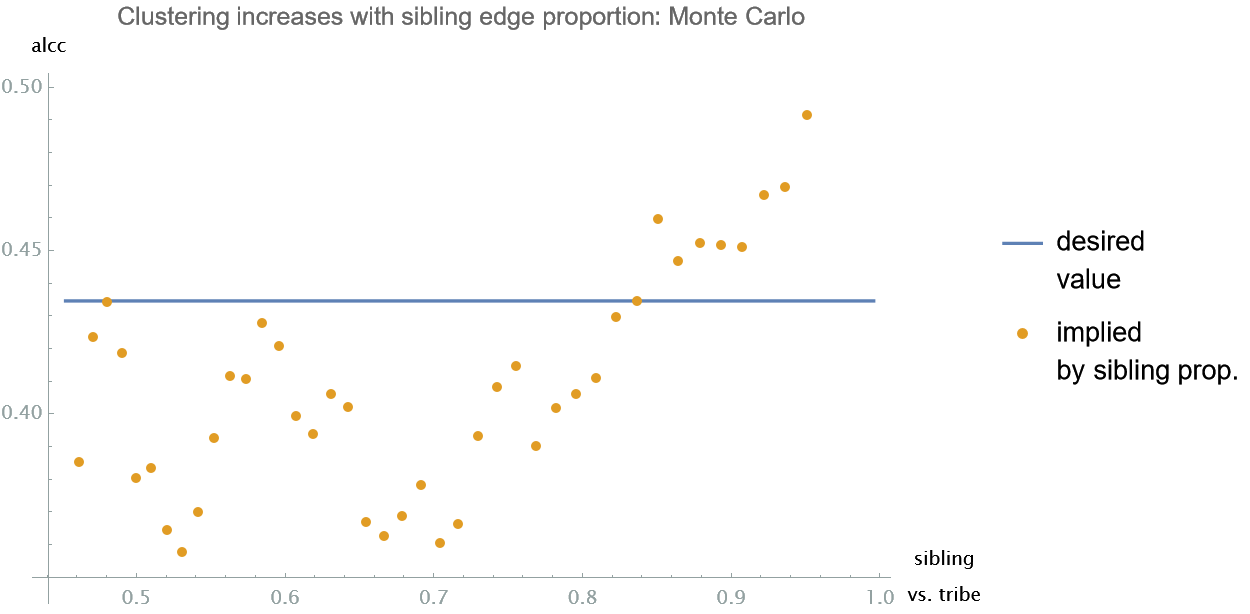}}

	\end{center}	
\end{figure}

\subsubsection{Mathematica script for \texttt{alccPredictor}} \label{sc:alccp}
\begin{Verbatim}[frame=single,fontsize=\small]
alccPredictor[Phi1_, Phi2_,  b1_, b2_, Sigma_] := 
   Module[{Phi, ePhi1, ePhi2, ePhi, q1, q2, tribalVertices, tribalMarks, 
                             edges1, edges2, tribalEdges},  
  tribalVertices = Array[v, {b1, b2}];   
  {ePhi1, ePhi2} =   {Exp[Phi1], Exp[Phi2]}; 
  ePhi = ePhi1 + ePhi2; 
  q1 = ePhi1/ePhi; 
  q2 = 1.0 - q1; 
  Phi = Log[ePhi]; 
  tribalMarks = RandomVariate[LogNormalDistribution[Phi, Sigma], {b1, b2}]; 
  edges1 = edgeGenerator[Flatten[tribalVertices], Flatten[tribalMarks], q1]; 
  edges2 =  Flatten[Table[edgeGenerator[tribalVertices[[k]], 
     tribalMarks[[k]], q2], {k, 1, b1}]]; 
  tribalEdges = loopRemover[Join[edges1, edges2]]; 
  Return[N[MeanClusteringCoefficient[Graph[tribalEdges]]]]; 
];
\end{Verbatim}
\subsubsection{Parsing the script \ref{sc:alccp} }
When leaf marks are Log-normal$(\phi, \sigma)$ random variates,
and the probability vector for the graph is $(q_0, q_1, q_2)$,
create new parameters
\[
\phi_1:=\phi + \log{q_1}; \quad \phi_2:=\phi + \log{q_2}
\]
which refer to tribal edges and clan edges, respectively, as in Lemma \ref{l:tribalsplit}.
To predict the ALCC in
the agreement graph, forget the root edges 
(since almost  none of them are in the agreement graph)
simulate a tribal graph, and compute its ALCC\footnote{
It is regrettable
that no effective analytic estimate of the ALCC could be found.
}.
The Script \ref{sc:alccp} achieves this, taking $(\phi_1, \phi_2)$
as arguments, and computing a tribal subgraph of a 
hidden ancestor graph from which root edges have been omitted.

\subsection{Fitting the agreement graph}\label{s:agreegraphfit}
The scripts \ref{sc:tribsurv}, \ref{sc:maxalcc}, and \ref{sc:alccp} are now combined
to fit the agreement graph, within the context of a tribal subgraph.

\subsubsection{Mathematica script for \texttt{tribalGraphFit}}\label{sc:tribalGraphFit}
\begin{Verbatim}[frame=single,fontsize=\small]
tribalGraphFit[b0_,b1_,b2_,Sigma_, dA_,KappaA_,Delta_]:=Module[
          {Rho1, Rho2,Kappamax,Kappa,Kappaprev,Phi,Phi1,Phi2},
  {Rho1,Rho2}=tribalSurvivalRates[b0,b1,b2,Sigma];
  Phi2 = Log[dA/Rho2] - Sigma*Sigma/2.0; (* Initialization *) 
  Kappamax = Median[Table[maxPossibleClustering[b0,b1,b2,Sigma, dA] , {11}]];
  {Kappa,Kappaprev} = {Kappamax,Kappamax}; (* Initialization *)  
  If[KappaA>Kappamax ,Return[Indeterminate],
    While[(Max[Kappaprev,Kappa]>KappaA)\[And](Phi2>0.0),
      Phi2 = Phi2 - Delta;
      Phi1 =Log[(1.0/Rho1)*(dA*Exp[-Sigma*Sigma/2] - Rho2*Exp[Phi2])];
      Kappaprev =Kappa;
      Kappa = alccPredictor[Phi1, Phi2,  b1, b2, Sigma] ;
    ]; (* end of While *)
    Return[{Phi1, Phi2, ( Kappa+Kappaprev)/2.0}];
  ]; (* end of If *)
];
\end{Verbatim}
\subsubsection{Parsing the script \ref{sc:tribalGraphFit} }
Check first that the actual value of $\kappa_A$ is less than $\kappa_{\max}$ returned by
\texttt{maxPossibleClustering}, for which we take the median of 11 estimates for robustness.
Then decrement the loading on clan edges (represented by $\phi_2$) in steps of size $\delta$ until the ALCC estimate from 
\texttt{alccPredictor} falls below $\kappa_A$ on two consecutive steps. Visualize this process
with Figure \ref{f:clustering}. Here the tribal graphs have only 2000 vertices, so ALCC Monte Carlo estimates
are noisy. The parameters returned are $q_1 e^{\phi}$, $q_2 e^{\phi}$ (compare Lemma \ref{l:tribalsplit}),
and the average of the last two values of ALCC.

\subsection{Combine conflict and agreement fits}
We have now reached the bottom node of the dependency diagram in Figure \ref{f:dependencies}.

\subsubsection{Mathematica script for \texttt{hagParameterFit}} \label{sc:hagfit}
\begin{Verbatim}[frame=single,fontsize=\small]
hagParameterFit[b0_,b1_,b2_,Sigma_,dC_,pC_,dA_,KappaA_,Delta_]:=Module[
         {Omega, Phi0, Phi1, Phi2,Phi, q0, q1, q2,Kappa},
  {Omega, Phi0} = conflictFit[b0,b1,b2,Sigma,dC,pC,Delta];
  {Phi1, Phi2, Kappa}=tribalGraphFit[b0,b1,b2,Sigma, dA,KappaA,Delta];
  Phi = Log[Exp[Phi0] +Exp[Phi1]+Exp[Phi2]];
  {q0, q1, q2} = {Exp[Phi0-Phi],Exp[Phi1-Phi],Exp[Phi2-Phi]};
  Return[{Omega, Phi, q0, q1, q2}];
];
\end{Verbatim}
\subsubsection{Parsing the script \ref{sc:hagfit} }
Script \ref{sc:confit} returns not only $\omega$, but also
$q_0 e^{\phi}$, while  script \ref{sc:tribalGraphFit} returns $q_1 e^{\phi}$ and $q_2 e^{\phi}$.
Since $q_0 + q_1 + q_2 = 1$, addition gives us $e^{\phi}$, and ratios give the values of 
$q_0, q_1, q_2$. These operations are performed by script \ref{sc:hagfit},
which completes the fitting of parameters $(\omega, \phi,q_0, q_1, q_2)$,
given observables $(d_C, p_C, d_A, \kappa_A)$, as sought in Definition \ref{d:identify}.

\subsection{Conclusion about simplified parameter identification} \label{s:solveidentify}
Script \ref{sc:hagfit} answers the question posed in Section \ref{s:simpleidentify}.\\
\textit{Script \ref{sc:hagfit}, and its dependencies as shown
in Figure \ref{f:dependencies}, together solve the simplified parameter identification problem
(Definition \ref{d:identify}) for hidden ancestor graphs.}

\section{Hidden ancestor graph parameter identification: general case}\label{s1:haggen}
\subsection{Modeling context}
Data engineers who gather intelligence from large labelled graphs typically deal with
data streams from which a new graph instance is built from a daily (or other periodic) batch of data.
These successive graph instances may share the same system of vertex labels, and may have vertices in
common, but the edges are new each time. Such instances may be appropriately modelled by 
hidden ancestor graphs with the same latent tree parameters $(b_0, b_1, b_2)$ and Log-normal leaf marks
with the same $\sigma$; if these are known, then solving
the parameter identification problem of Definition \ref{d:identify}, as we described 
in Section \ref{s:solveidentify} may suffice.

This section describes how to make an initial choice of 
latent tree parameters $(b_0, b_1, b_2)$ and Log-normal parameter $\sigma$, in
the context of the Poisson Log-normal model used for hidden ancestor graph generation.

\subsection{Missing vertices in the agreement graph}
We describe in Section \ref{s:quartilefit} of the Appendix a method of fitting a Poisson Log-normal 
model to the agreement graph of the real graph instance. Suppose the fitted parameters are
$(\phi_A, \sigma)$. The $\phi_A$ parameter cannot be used directly, 
since the edge losses due to loops are as yet unknown. However these edge losses act as
a thinning parameter on Poisson variates with Log-normal distribution,
and the $\sigma$ parameter is not affected by this thinning. Hence we may use it in the scripts of Section
\ref{s1;hagparamid}.

This fitting also estimates the number of invisible vertices,
i.e. those which did not acquire any half-edges because their
leaf mark was too small. See Script \ref{sc:logRatioError}.
By adding this number to the count of visible vertices, we obtain a suitable
number $n$ for the vertex count $n:=b_0 b_1 b_2$ of a hidden ancestor graph.

It is crucial in fitting a hidden ancestor graph model to estimate the proportion $\alpha$ of missing vertices
-- typically as large as 25\% -- which is needed to inflate the denominator when estimating
the mean agreement degree $d_A$ and mean conflict degree $d_C$. For example, suppose 21 million agreement edges
are observed among 0.75 million occupied vertices in a real graph where we estimate $\alpha = 0.25$. The average agreement degree
per vertex (including missing vertices) is not $(2 \times 21) \div 0.75 = 56.0$, but $42 \div 1.0 = 42.0$, with a similar adjustment for conflict degree.
Failure to estimate $\alpha$ also leads to overestimates of Log-normal parameters.

\subsection{Latent tree parameters}
\subsubsection{Number of tribes}
It is natural to take the number $b_0$ of tribes in the latent tree
to be the number of distinct labels in the real graph instance.

\subsubsection{Number of siblings in a clan}
Subject matter expertise may suggest a suitable value for $b_2$ 
for a given context. If $b_2$ is chosen too large,
then script \ref{sc:maxalcc} will return a maximum possible ALCC less than the $\kappa_A$
observed in the real graph. A small value of $b_2$
may lead to a model which can be fitted, but whose generation is inefficient because of an abundance of loops.
Simulation experiments when the \texttt{edgeGenerator} was applied
to $n$ Log-normal$(2.0, 2.1)$ leaf marks showed that as $n$ increased
from 50 to 2500, average local clustering in the visible graph
decreased from about 0.63 to about 0.5, showing that small clan
sizes are not essential for high average local clustering.

\subsubsection{Number of clans in a tribe}
Now that we can approximate $n:=b_0 b_1 b_2$, and we know $b_0$ and $b_2$, setting
$b_1:=\lfloor[n/(b_0 b_2)\rfloor$ yields a hidden ancestor graph of about the same size
as the real graph instance. On the other hand, if the real graph instance has $10^7$ vertices,
and we would rather test algorithms on graphs with $10^5$ vertices, then choose $b_1$ 100 times smaller.
This will have a modest effect on observable graph parameters, other than adding more noise.

\subsection{Conclusion}
The latent tree parameters $(b_0, b_1, b_2)$ and Log-normal $\sigma$ can be estimated using
methods described here, as a precursor to those of Section \ref{s1;hagparamid}.

\section{Case Study}\label{s:casestudy}
\subsection{Citation data}
The purpose of this section is to demonstrate fitting of a hidden ancestor graph model to a real world graph constructed from a citation network.
For this study we used the DBLP Citation Network v12 dataset \cite{dblp}; here nodes of the graph are journal articles, and
edges are based on co-author links, as in Table \ref{tab:divers}.
We extracted the agreement graph using the edges where end points were papers on the same ``topic",
and the conflict graph using the edges where the ``topics" were different; here topic refers to one of 9 fields of study 
from the database. We computed the proportion of vertices appearing in the conflict graph, the average clustering coefficient of the agreement graph,
the degree distribution and its quantiles, and the Log-normal fit parameters. These statistical estimates are shown in Table~\ref{tab:dblp_params}.
This analysis does not adjust for invisible vertices, because the later simulation showed there are few of them.

\begin{table}[htbp]
    \centering
    \begin{tabular}{c|c}
         Number of nodes & 4,490,992\\
         Number of edges & 295,035,281\\
         Vertex degree Log-Normal $\sigma$ & 1.60\\     
         Avg agreement \& conflict degrees $(d_A, d_C)$ & $(101.081, 30.1325)$ \\
         Proportion in conflict graph $p_C$ & 0.89\\
         Avg clustering (agreement graph)  $\kappa_A$ & 0.762\\ 
    \end{tabular}
    \caption{\textit{Summary statistics from DBLP V12 Co-author graph.}}
    \label{tab:dblp_params}
\end{table}
Using the code described in Sections \ref{s:conflgraphfit}  and \ref{s:agreegraphfit} we 
fitted a 1\% scale hidden ancestor graph model to these statistics. Model parameters are shown in
Table \ref{tab:dblp-hag}. Two repetitions led to nearly identical parameters. We took $b_0 = 9$ as the number of topics (i.e. tribes),
and $b_2 = 25$ because larger values did not allow $\kappa_A$ as big as $0.76$. 
\begin{table}[htbp]
    \centering
    \begin{tabular}{c|c}
      Latent tree parameters $(b_0, b_1, b_2)$  &  $(9, 220, 25)$\\
      Wildness rate $\omega$  & $0.201$ \\
      Log-normal parameters $(\phi, \sigma)$ & $(4.14, 1.6)$ \\
      Probability vector $(q_0, q_1, q_2)$ & $(0.417, 0.003, 0.581)$\\
    \end{tabular}
    \caption{\textit{Fitted hidden ancestor graph parameters for DBLP V12 Co-author Model}}
    \label{tab:dblp-hag}
\end{table}

Finally we generated edges for this graph on $9 * 220*25 =49500$ vertices, and measured the graph parameters.
Observe that 99\% of all vertices actually appear in the agreement graph.
Results are satisfactory except that the average clustering coefficient (0.48) in the agreement graph fell below the desired target.
Possibly our notion of ``clan'' is too crude as a model of scientific co-authorship.
\begin{table}[htbp]
    \centering
    \begin{tabular}{c|c}
          Total multi-edges (excl. loops) &  3,444,719 \\
      Agreement graph & 48,789 vertices, 527,247 unique edges  \\
Conflict graph & 43,347 vertices, 683,873 unique edges  \\
Avg agreement \& conflict degrees $(d_A, d_C)$ & $(109.8, 27.63)$ \\
Proportion in conflict graph $p_C$ & 0.88\\
Avg clustering (agreement graph)  $\kappa_A$ & 0.48\\ 
    \end{tabular}
    \caption{\textit{Hidden ancestor graph 1\% scale simulation of DBLP V12 Co-author graph}}
    \label{tab:dblp-simulate}
\end{table}

\section{Binary classifier for wild vertices}\label{s1:class}
\subsection{Alternative approaches for building a classifier}
The purpose of this section is to build a binary classifier of the hidden state of each conflict vertex
-- wild or tame -- from the edges and vertex labels of the graph.
\subsubsection{Structural approach}
Our approach will be to infer structural parameters of a hidden ancestor graph model from the
edge and label data, and then to create a predictor based on the inferred model.
The specificity of the hidden ancestor graph model of Section \ref{s1:hag} permits
a model combining a natural Bayes formula and a simple combinatorial principle.
It involves a logistic predictor, but this predictor is unsupervised: it need {\em not} be trained
on test data but has coefficients derived from fitted parameters of the hidden ancestor graph.
First we contrast this approach with a couple of model-agnostic alternatives.
\subsubsection{Supervised learning}
A statistician tasked with building a binary classifier will often 
train a supervised learning model. Start by generating
 two instances of hidden ancestor graphs with the same underlying parameters.
In one of these graphs, wild status is revealed, and is used to train a binary classifier,
which is evaluated on the other graph. Since the domain consists of graph vertices, a
statistician may attempt a graph neural network model \cite{err}.
In practice, colleagues who implemented a graph neural network classifier encountered
much heavier computing loads and inferior performance to a Bayes belief propagation model.
\subsubsection{Bayes belief propagation}
Attach to every vertex and both ends of every edge a log Bayes ratio,
which is updated via message passing between vertices and proximal ends, and between
the two ends of an edge. In a sequence of message passing rounds,
each vertex's belief in its probability of wild status is updated.
Agreement edges have an effect of pushing this belief towards zero, while conflict edges
push this belief towards one in a vertex whose conflict neighbors believe they are tame.
There is no theoretical convergence guarantee, but in practice message passing converges
in twenty rounds or less.
Bayes belief propagation has the advantage of being adaptable to
inhomogeneity in wildness rates among different tribes, for example.

\subsection{Vertex covers}
We repeat that our goal here is not community detection (unnecessary because tribal marks
are visible) but the binary classification problem of deciding which vertices in the conflict graph
are wild. By construction of the conflict graph, where every conflict edge has at least one wild end point,
the wild vertices constitute a vertex cover for the conflict edges. However they are not usually the smallest
vertex cover. 

A reasonable approach is to choose a vertex cover $U$ such that the expected number of tame vertices in $U$ is as small
as possible. One approach is to derive a \textit{vertex weight function} for conflict vertices that
estimates the conditional probability that a vertex is tame, given observable statistics such as the
conflict valency and agreement valency of the vertex.
If we find a vertex cover of minimum weight with respect to such a function,
classifying every vertex in this cover as wild minimizes the expected number of mis-classifed tame vertices.

\begin{figure}
    \caption{\textit{Associate to each vertex in the conflict graph the
    pair consisting of agreement valency ($x$-axis) and conflict valency ($y$-axis).
    Wild and tame vertices supply separate scatter plots.
    The upper and lower curves shows $(x,y)$ pairs for which the tameness probability
    (\ref{e:conditionaltamecompute}) takes values 0.01 and 0.99, respectively. Their placement
    suggests the power of the Logistic predictor to separate wild and tame vertices of higher valency.
    }
    }
    \label{f:wildtamesplit}
    	\begin{center}
\scalebox{0.45}
{\includegraphics{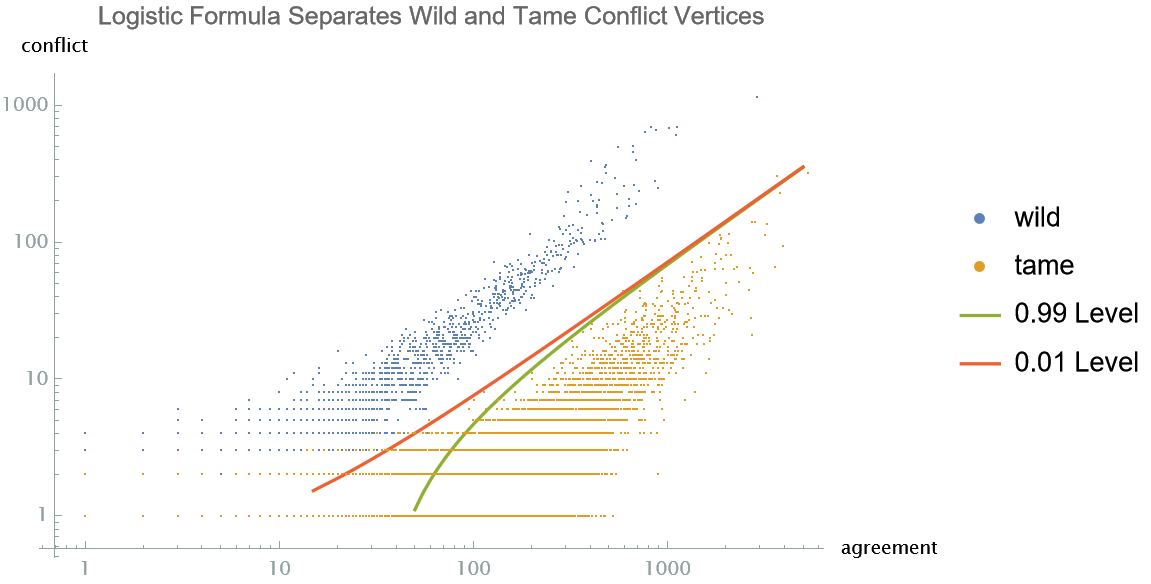}}
	\end{center}	
\end{figure}

A natural vertex weight function heuristic is described in Section \ref{s:weightfunction}.
Figure \ref{f:wildtamesplit} demonstrates by example that the bivariate random vector consisting of
agreement valency ($x$-axis) and conflict valency ($y$-axis) indeed has a distribution which depends
significantly on whether the vertex is wild or tame.

\subsection{Tameness probability: Logistic predictor}\label{s:weightfunction}
Suppose a vertex of the hidden ancestor graph has $X_A$ incident agreement edges,
and $X_C$ incident conflict edges.
Now we propose a specific logistic predictor\footnote{
{\em Not} trained
on test data: coefficients are derived from fitted parameters of the hidden ancestor graph.
} of tameness probability.
Figure \ref{f:wildtamesplit} plots a level curve of this formula, and
lends empirical support to its power to separate wild and tame vertices.

\subsubsection{Logistic predictor: }

\textit{Approximate the conditional probability
that a conflict vertex is tame by:}
\begin{equation}\label{e:conditionaltamecompute}
    \pr[\mbox{tame} \mid X_A = a, X_C = c] \approx
    \frac{1-\omega}
    {1-\omega + \omega \exp({- a (1-\omega) \lambda - c \log{\omega})}}.
\end{equation}
\textit{where  $\lambda: = q_0 / (q_1 + q_2)$}.

\begin{figure}
    \caption{
    \textit{ Logistic predictor (\ref{e:conditionaltamecompute}) was calculated
    for 2K tame and 2K wild vertices, and the empirical CDF was computed for each.
The orange curve shows that tameness probability ($x$-axis)
was predicted to be less than 0.1 ($y$-axis) for nearly 90\% of wild vertices.
The blue curve shows it to be at least 0.5 for about 70\% 
of tame vertices. This accords with the wild vs. tame
separation seen in Figure \ref{f:wildtamesplit}.
    }
    }
    \label{f:histogramstameprob}
    	\begin{center}
\scalebox{0.5}
{\includegraphics{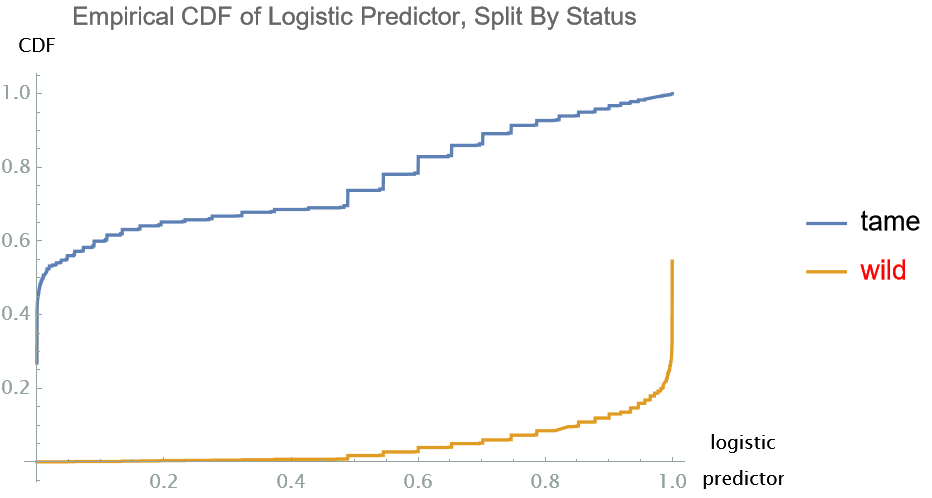}}
	\end{center}	
\end{figure}

This formula is valid even when $a=0$ (a conflict vertex with no incident agreement edges).
Figure \ref{f:histogramstameprob} shows how well (\ref{e:conditionaltamecompute}) 
discriminates in practice between wild and tame vertices.
The right side of (\ref{e:conditionaltamecompute}) is an algebraic simplification of the
expression:
\begin{equation}\label{e:conditionaltame}
    \frac{(1-\omega) P(c, \omega \lambda a)}
    {\omega P(c, \lambda a)  + (1-\omega) P(c, \omega \lambda a)},
    \quad a \geq 1, c \geq 1.
\end{equation}
where $P(k, \theta)$ refers to the Poisson probability $e^{-\theta} \theta^k / k!$.
The rationale for (\ref{e:conditionaltame}) arises from Lemma \ref{l:tamegivenmark}, 
which is based on the fundamental hidden ancestor graph construction.

\begin{lemma}\label{l:tamegivenmark}
Neglecting loops and agreement edges among root edges,
the conditional probability that $v$ is tame, given $X_C = c$ conflict edges
incident to $v$, and mark $F$ on $v$, is
\begin{equation}\label{e:tamegivenmark}
    \frac{(1 - \omega)  P(c, \omega q_0 F)}{(1 - \omega)  P(c, \omega q_0 F) + \omega P(c, q_0 F)}.
\end{equation}
\end{lemma}

\begin{proof}
Condition on the event that vertex $v$ has leaf mark $F$. The number of root half-edges associated with
$v$ has a Poisson$(q_0 F)$ distribution. Neglecting loops and agreement edges among root edges,
the conflict valency of $v$ is Poisson$(q_0 F)$ if $v$ is wild (because half-edges associated with
$v$ can be paired with those of any tribe), but Poisson$(\omega q_0 F)$ if $v$ is tame, because
pairings are possibly only when the other vertex is wild, which has probability $\omega$.
In other words:
\[
\pr[X_C = c \mid F, \mbox{wild}] = P(c, q_0 F); \quad 
\pr[X_C = c \mid F, \mbox{tame}] = P(c, \omega q_0 F].
\]
Use Bayes' Rule:
\[
\pr[\mbox{tame} \mid F, X_C = c] = \frac{\pr[\mbox{tame} \cap \{X_C = c\} \mid F]}{\pr[ \{X_C = c\} \mid F]}
= \frac{(1 - \omega)  P(c, \omega q_0 F)}{\pr[ X_C = c \mid F]}.
\]
The denominator is $(1 - \omega)  P(c, \omega q_0 F) + \omega P(c, q_0 F)$ by conditioning on whether
$v$ is wild or tame.
\end{proof}

\subsubsection{Justification for the heuristic (\ref{e:conditionaltame})}
\label{s:justify}
\begin{figure}
    \caption{\textit{
    A scatter plot supports the heuristic of Section \ref{s:justify}.
    Agreement valence, boosted by $1/ (q_1 + q_2)$, is shown
    on the $x$-axis, and leaf mark on the $y$-axis.
    The relationship is approximately linear (with slope of 1),
    but weaker at the low end because of variance, and at the high end
    because of loop removal for high leaf marks.
    } }
    \label{f:markratio}
    	\begin{center}
\scalebox{0.5}
{\includegraphics{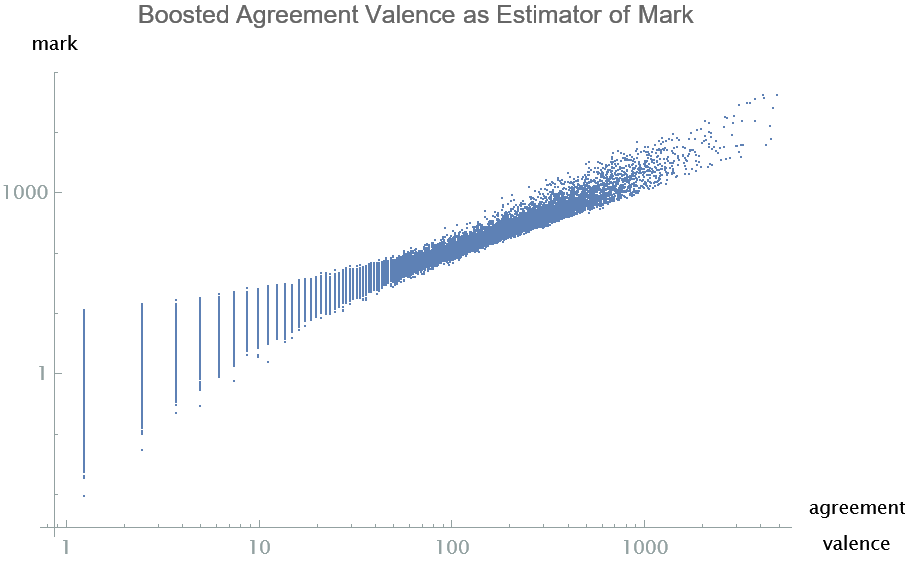}}
	\end{center}	
\end{figure}

We cannot use (\ref{e:tamegivenmark}) directly to predict the probability of tameness, because
the leaf mark $F$ is unobserved. However the number $X_A$ of agreement edges incident to $v$ is independent
of whether $v$ is wild or tame, and is roughly Poisson $((q_1 + q_2)F)$, neglecting loss due to loops,
and neglecting root level agreement edges. The heuristic is to replace $q_0 F$ by $q_0 X_A / (q_1 + q_2)=: \lambda X_A$.
Here $X_A$ is observed, and $(q_0, q_1, q_2)$ and $\omega$ may be estimated as above.
The validity of such a substitution (or lack of it) can be assessed from 
Figure \ref{f:markratio}, where $X_A/(q_1 + q_2)$ is shown on the $x$-axis, and $F$ on the $y$-axis, in log-log scale.
There are two main sources of error in the substitution.
\begin{enumerate}
    \item When $F$ is small, $X_A$ has variance close to $(q_1 + q_2)F$.
    \item When $F$ is large, $X_A / F$ is reduced by loops at the sibling level.
\end{enumerate}
Despite these sources of error, 
the effectiveness of the heuristic is evident in the scatter plot of Figure \ref{f:wildtamesplit}
and the cumulative distributions in Figure \ref{f:histogramstameprob}.

\subsubsection{Other ways to fit a logistic predictor}
The coefficients $\omega$ and $\lambda: = q_0 / (q_1 + q_2)$ in (\ref{e:conditionaltame}) can be fitted by the methods of Sections
\ref{s1;hagparamid} and\ref{s1:haggen},  but this is not the only way to obtain a logistic predictor of tameness.
An alternative is to fit a hidden ancestor graph model to graph data,
generate an instance $H$ of model, sample equal numbers of wild and tame
vertices in $H$, and fit a logistic regression in the usual way
to obtain parameters $\beta_0, \beta_1, \beta_2$ in a formula:
\[
\frac{1}{1 +  \exp{(-\beta_0 - \beta_1 a - \beta_2 c)}}.
\]

\subsection{Rapid binary classifier of wild status}\label{s:rapidclassify}
Consider conflict graphs arising from parameters in Table \ref{tab:genericparameters}.
Suppose the statistician wishes to determine the wild status of conflict vertices.
The classifier (\ref{e:conditionaltamecompute}) cannot be used in isolation,
because it may classify both end points of a conflict edge as tame, which
violates the rules of hidden ancestor graphs. However (\ref{e:conditionaltamecompute}) 
can serve as vertex weighting for the {\em minimum weight vertex cover}
problem, treated fully in standard texts such as Vazirani \cite{vaz}. Typical greedy algorithms
make a single pass through the conflict vertices, making updates to
the conflict graph at every step.

Since (\ref{e:conditionaltamecompute}) is a  probability,
values very close to zero or one alone contain strong classification information.
Indeed for the examples shown in Table \ref{tab:rapidclassify}, tameness probability (\ref{e:conditionaltamecompute}) 
was at least 0.99 for 67\% of conflict vertices,
and under 0.01 for 8\% of conflict vertices.
This leads to a {\em rapid classifier} much faster and simpler than greedy
one-vertex-at-a-time algorithms.

Suppose that a conflict graph is supplied, as part of a hidden ancestor graph
for which parameters $(q_0, q_1, q_2, \omega)$ have been estimated, so
that the Logistic predictor (\ref{e:conditionaltamecompute}) can be evaluated
for every conflict vertex. The following steps are mirrored in Table \ref{tab:rapidclassify},
and recapitulated in Script \ref{sc:rapidclassify}.

\begin{table}[]
    \centering
    \begin{tabular}{c||c|c|c|c|c}
    \toprule
Classification & WHP  & Adj. to Tame & Singletons & Not classified & Mis-classified \\ \midrule
 Wild  &  1,173 & 1,102 & \, & (62) & 51\\
Tame & 12,365 & \, & 3,732 & (64) & 8 \\ \bottomrule
    \end{tabular}
    \caption{Hidden ancestor graph parameters come from Table \ref{tab:genericparameters}.
    Out of 18,642 conflict vertices, all but 126 received
    a classification using the rapid classifier, and all but 59 of the remaining 18,514 
    were classified correctly. The WHP column reflects initial judgement
    where (\ref{e:conditionaltamecompute}) is less than 0.01 or
    greater than 0.99. The next column shows vertices marked as wild by
    adjacency to a WHP tame vertex. In the subgraph on remaining vertices,
    singletons were declared tame, and 126 vertices without a classification 
    were in components
    of size 2 or 3. Figure \ref{f:rapidclassify} illustrates similar data.}
    \label{tab:rapidclassify}
\end{table}

\begin{enumerate}
    \item[(a)] Select a small $\epsilon$, such as 0.01.
    \item[(b)] Mark as wild all vertices for which (\ref{e:conditionaltamecompute}) is less than $\epsilon$.
    \item[(c)] Of the vertices for which (\ref{e:conditionaltamecompute}) exceeds $1-\epsilon$, select
    a maximal independent set\footnote{
    Identify the handful of edges connecting vertices with high tameness probability, compute a minimal vertex cover, and
    exclude vertices in the cover from the tame list.
    } in the conflict graph, and declare all its elements tame.
    \item[(d)] Classify as wild any
     remaining vertex adjacent to any already marked as tame.
    \item[(e)] In the subgraph $H$ induced by vertices not already classified,
    singletons are declared tame (since they have no incident conflict edges with
    a tame vertex at the other end). Others are ``not classified''.
\end{enumerate}
As we see in Table \ref{tab:rapidclassify}, the subgraph $H$ also contains
components of size 2 or more. These may be small enough that exhaustive minimum
weight vertex cover algorithms may be applied, but for now we shall leave them as ``not
classified''. Larger values of $\epsilon$ reduce
the number of these ``not classified'' vertices, while raising the number of mis-classified ones by a similar amount.

\subsubsection{Mathematica script for rapid classifier}\label{sc:rapidclassify}
\begin{Verbatim}[frame=single,fontsize=\small]
wildTameUnknown[cgraph_, q0_, q1_, q2_, omega_, eps_]:=Module[{tameWHP, wildWHP,
excluded, tameIndep, middle, wildByAdj, leftovers, components, tameSingletons}, 
    tameWHP = Select[VertexList[cgraph], 
      tamePr[degA[#1], degC[#1], q0, q1, q2, omega] > 1. - eps & ]; 
    excluded = FindVertexCover[Graph[EdgeList[Subgraph[cgraph, tameWHP]]]]; 
    tameIndep = Complement[tameWHP, excluded]; 
    wildWHP = Select[VertexList[cgraph], 
      tamePr[degA[#1], degC[#1], q0, q1, q2, omega] < eps & ]; 
    middle = Complement[VertexList[cgraph], tameIndep, wildWHP]; 
    wildByAdj = Select[middle, 
      Length[Intersection[AdjacencyList[cgraph, #1], tameIndep]] > 0 & ]; 
    leftovers = Complement[middle, wildByAdj]; 
    components = ConnectedComponents[Subgraph[cgraph, leftovers]]; 
    tameSingletons = Flatten[Select[components, Length[#1] == 1 & ]]; 
    Return[{{wildWHP, wildByAdj}, {tameIndep, tameSingletons}, 
      Complement[leftovers, tameSingletons]}]; 
 ]; 
\end{Verbatim}
\subsubsection{Parsing the script \ref{sc:rapidclassify}}
\texttt{degA} and \texttt{degC} are functions returning the agreement and conflict valence,
respectively, while \texttt{tamePr} refers to (\ref{e:conditionaltamecompute}).
The script returns wild vertices of the types (b) and (d), tame vertices of the types (c) and (e), 
and ``not classified'' vertices,
respectively, as seen in Table \ref{tab:rapidclassify}.
The union of \texttt{tameIndep} and \texttt{tameSingletons}
forms an independent set in the conflict graph, though not necessarily a maximal one.
\subsubsection{Performance of the rapid classifier}
\begin{figure}
    \caption{\textit{Rapid classifier based on logistic predictor:
    relative sizes of the five sorts of outcomes for conflict vertices,
    as in Table \ref{tab:rapidclassify}
    }
    }
    \label{f:rapidclassify}
    	\begin{center}
\scalebox{0.5}
{\includegraphics{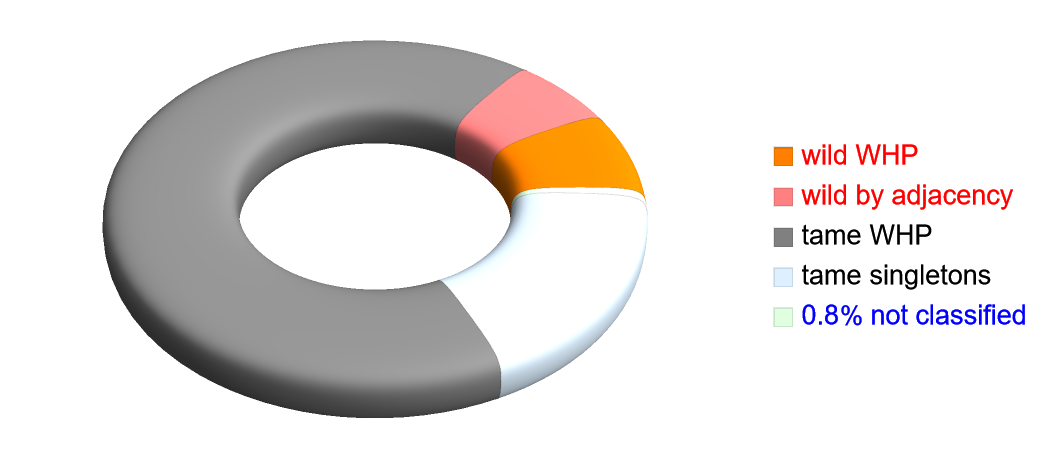}}
	\end{center}	
\end{figure}

Table \ref{tab:rapidclassify} shows that, for hidden ancestor graph examples
of the type of Table \ref{tab:genericparameters}, the rapid classifier with $\epsilon = 0.01$
reached a verdict for over $99\%$ of conflict vertices, with error
rates of about 2\% on wild vertices, and about 0.2\% on tame ones. 
Here is the confusion matrix for an example with $b_1 = 100$ (instead of 50),
giving a graph with 35,559 conflict vertices, 35,016 of which are in a giant component.
\begin{equation}\label{e:confusionmat}
  \begin{matrix}
\mbox{conflict graph} & \mbox{classified wild} & \mbox{classified tame}& \mbox{not classified} & \mbox{TOTAL}\\ \hline
 \mbox{true wild} & 4902 & 81 & 139 & 5122 \\
  \mbox{true tame}  & 14  & 30320 & 143 & 30477 \\ \hline
  \mbox{totals}  & 4916 & 30401 & 282 & 35559
\end{matrix}
\end{equation}
In this example:
\begin{enumerate}
\item 
The precision of the wild classifier is $\frac{4902}{4916}=99.7\%$ and the recall is $\frac{4902}{5122}=95.7\%$.
    \item
    Of 4,916 vertices classified as wild, 2,623 come from (b) and 2,293  from (d).
\item    
Of the 30,401 classified as tame, 23,699 come from (c) and 6,702 from (e).
\end{enumerate}
Similar statistical ratios are seen in larger examples. The pie chart in Figure \ref{f:rapidclassify}
shows typical relative sizes of the five sorts of outcomes for conflict vertices.

\subsubsection{How rapid classifier performs on conflict graph 2-core}
In the example shown in Table \ref{tab:rapidclassify}, 45\% of conflict vertices are
in the 2-core of the conflict graph. Apply the rapid classifier to the entire conflict graph,
but focus on its verdict on vertices in the 2-core of the conflict graph:
only 10 vertices in the 2-core were not classified, and
the confusion matrix has even lower false positive and false negative rates on the 2-core
shown in (\ref{e:confusionmat2core}):
\begin{equation}\label{e:confusionmat2core}
  \begin{matrix}
\mbox{2-core only} & \mbox{classified wild} & \mbox{classified tame}& \mbox{not classified}& \mbox{TOTAL}\\ \hline
 \mbox{true wild} & 1719 & 0 & 2 & 1721\\
  \mbox{true tame}  & 7  &  6722 & 8 & 6737\\ \hline
  \mbox{totals}  & 1726 & 6722 & 10 &8458
\end{matrix}  .
\end{equation}
Within the 2-core precision of the wild classifier is about the same, 
but the recall rises to $\frac{1719}{1721}=99.9\%$.
This example implies that most of the ``not classified'' and wrongly classified vertices are in the complement
of the conflict graph 2-core.

\section{Topics for future study}

\subsection{Demonstrate consistency of parameter estimates}\label{s:consistency}

The justification for using observables $(d_C, p_C, d_A, \kappa_A)$ to estimate
$(\omega, \phi, q_0, q_1, q_2)$ for known $(b_0, b_1, b_2,\sigma)$
could be rigorously established by proving
the following asymptotic conjecture. The methods of Bickel and
Chen \cite{bic} may be applicable.

\textbf{Conjecture: }\textit{
Suppose marks are i.i.d. Log-normal$(\phi, \sigma)$ random variates. Let $b_1 \to \infty$ holding 
$(b_0, b_2,\omega, \phi, \sigma, q_0, q_1, q_2)$ fixed. Then the statistics $(d_C, p_C, d_A, \kappa_A)$ of
the resulting sequence of hidden ancestor graphs converge in law to a limit in $R^4$ which, for each $b_0, b_2$ and $\sigma$, is uniquely
determined by parameters $(\omega, \phi, q_0, q_1, q_2)$ (with constraint $q_0 + q_1 +q_2 = 1$).
}

\subsection{Inhomogeneous structure in conflict graph}
For simplicity the focus has been on depth three hidden ancestor graphs, which have the
merit of having identifiable parameters. So far as realism is concerned, they have a major drawback.
Consider the edge-weighted {\em supergraph} whose nodes are the $b_0$ tribal labels, and where tribes $i$
and $j$ are linked by an edge whose weight $w_{i,j}$ is the sum of all conflict graph edges between
vertices in tribes $i$ and $j$ respectively. For large values of $b_1 b_2$ this will look like a complete weighted graph
with similar weights on all the edges. In vertex labelled graphs 
from  applications described in Table
\ref{tab:divers} it is likely that certain pairs of tribes will show much more conflict than
other pairs. This may reflect hierarchy in the partitioning into tribes, as can be seen from the examples of Table \ref{tab:divers}:

\textit{Genomics: }Phenotypes fall into a hierarchy based on the functional components (respiration,
perception, digestion, mobility, etc.) of the organism. 

\textit{Retailing: }The label attached to each item is a consumer category, which has a hierarchical structure, such as
\textit{DIY hardware $\rightarrow$ hammer}, or \textit{apparel $\rightarrow$ mens $\rightarrow$ boots}. 

\textit{Telecoms: }Service areas follow a geographic hierarchy such as \textit{country $\rightarrow$ province $\rightarrow$
city $\rightarrow$ postcode}.

\textit{Bibliography: }
The ACM Computing Classification System,
has a hierarchical structure such as D $\rightarrow$ D.2 $\rightarrow$ D.2.2.

An enhanced model would create more realistic supergraph structure. There are at least two variations to consider.

\subsubsection{Unequal tribe sizes}
In the model of Section \ref{s1:hag}, every tribe has exactly $b_1$
clans. Instead a probabilistic partitioning scheme such as the
Chinese Restaurant Process \cite{pit} could be used: treat the
$n$ vertices as ``customers'' who enter a restaurant, and are
assigned in groups of $b_2$ to ``tables''. Stop when $b_0$ tables
(i.e. tribes) are occupied. This gives an approximately Dirichlet 
distribution of clans per tribe. Methods of Section \ref{s1;hagparamid}
no longer provide parameter identification. 

\subsubsection{Models of depth four or more}\label{s:hagdepth4}
A hidden ancestor graph of depth four would depend on integers $(b_0, b_1, b_2, b_3)$,
where tribes (whose labels are visible on vertices) correspond to $b_0 b_1$ nodes at depth two.
The ancestor has $b_0$ chiefdoms, and each chiefdom has $b_1$ tribes. There are
four levels of edge generation, instead of three. The tribal and clan edges are agreement edges as before,
but root edges and chiefdom edges are mostly conflict edges. A chiefdom conflict edge will have end points
in different tribes, but the same chiefdom. Then the supergraph of the conflict graph will look
like $b_0$ clusters with $b_1$ elements in each.

\section {Conclusions}

\subsection{Where the hidden ancestor graph fits in the modelling portfolio}
Hidden ancestor (multi)graphs provide models where vertices are labelled
by their community membership (tribe), but possess a hidden state
of tribal disaffiliation. An edge between members of different tribes
is possible only when at least one of these members has disaffiliated.

\subsection{Identifiability of parameters in fitting the model}
These models can be generated at billion-edge scale, and can
be fitted to have Log-normal vertex degree distribution, with
high average local clustering and with a conflict subgraph of
the desired size.

\subsection{Rapid classification of wild and tame vertices using Logistic predictor}
The statistical design of hidden ancestor graphs permits a
logistic function of agreement degree and conflict degree at a vertex
to predict the probability of tribal disaffiliation.
Suitable coefficients for the logistic function can be inferred
from model parameters. A straightforward four-step 
combination of statistical and combinatorial operations
allows nearly all the disaffiliated vertices to be identified
with high precision.

 \subsection{Use cases}
Besides the case study in Section \ref{s:casestudy}, Table \ref{tab:divers} suggests that 
hidden ancestor graphs have the potential to model labelling anomaly problems in retailing, telecoms, informatics, and genomics.

\subsection*{Acknowledgments}
The authors received valuable suggestions from Tim McCollam, Emma Cohen, and Jon Berry; and corrections from
Tim Chow, Doug Jungreis, and Kevin Chen.

\appendix

\section{Fitting the zero-censored Poisson Log-normal model}\label{a:quartilefit}

\subsection{The censoring problem}\label{a:censor}
Suppose $Y$ is standard normal, so that the random variate $e^{\mu + \sigma Y}$ is 
Log-normal$(\mu, \sigma)$; then sample a Poisson random variable $Z$ with conditional mean $e^{\mu + \sigma Y}$.
Call $Z\in \{0, 1, 2, \ldots\}$ a Poisson Log-normal random variable.
Repeat this construction for i.i.d. Normal$(\mu, \sigma)$ $Y_1, Y_2, \ldots, Y_n$ times, 
without recording the values $(Y_i)$, and sampling a new Poisson variate $Z_i$ for each, but
discarding it if $Z_i=0$. At the conclusion the value of $n$ is forgotten, and
only strictly positive $(Z_i)$ remain. Such a situation occurs when a graph is generated whose 
vertex degrees follow a Poisson Log-normal model,
but isolated vertices are not counted.

Let $\alpha$ denote the unknown proportion of the Poisson variates discarded. For example, if
$(\mu, \sigma)=(2.6, 5.6)$, then $\alpha \approx 0.293$. 
There is no closed form expression for $\alpha$ in terms of $\mu$ and $\sigma$,
since the moment generating function of the Log-normal is intractable \cite{asm}.
Our goal is to estimate $(\mu, \sigma)$ from the collection of strictly positive $(Z_i)$,
without knowing $\alpha$. In this case the standard maximum likelihood estimates based on means and standard deviations of the logs of the Poisson variates will have significant bias, because
the lower tail of $(Y_i)$ is missing from the sample, and we do not have the correct denominator $n$.

In his Appendix, Bulmer \cite{bul} sketches a maximum likelihood estimation of parameters $\mu$ and $\sigma$
based on a censored (i.e. omitting zero values) Poisson Log-normal sample. However
Serfling \cite{ser} points out the lack of robustness of such methods in the presence of contamination
by outliers.
Clauset, Shalizi, and Newman \cite{cla} present other methods for fitting the Log-normal distribution
to graph data.
We propose a compute-intensive numerical approach based on quartiles which explicitly estimates $\alpha$
via three short scripts as follows.

\subsection{Monte Carlo method for quartile-based fitting }\label{s:quartilefit}
For the censored Poisson Log-normal with $(\mu, \sigma)=(2.6, 5.6)$, the quartiles based on a 50K simulation were 
$(Q_1, Q_2, Q_3):=(9, 109, 2559)$. If we knew $\alpha$, 
we could use script \ref{sc:lognormalestim} to estimate $\mu$ and $\sigma$.
\subsubsection{Mathematica script for \texttt{Log-normalParams}}\label{sc:lognormalestim}
\begin{Verbatim}[frame=single,fontsize=\small]
LognormalParams[Q1_, Q2_, Q3_, Alpha_]:= Module[
    {x1, x2, x3, Mu, Sigma},
  x1 = InverseCDF[NormalDistribution[], 0.75*Alpha + 0.25];
  x2 = InverseCDF[NormalDistribution[], 0.5*Alpha + 0.5];
  x3 = InverseCDF[NormalDistribution[], 0.25*Alpha + 0.75];
  Sigma = Log[N[Q3/Q1]]/(x3 - x1);
  Mu = Log[N[Q2]] - Sigma*x2;
  Return[{Mu, Sigma}];
  ]
\end{Verbatim}
\subsubsection{Parsing the script \ref{sc:lognormalestim} }
The arguments are the censored quartiles $Q_1, Q_2, Q_3$ and the censoring rate $\alpha$.
The values $x_1, x_2, x_3$ are quantiles of $N(0,1)$ corresponding to
$\alpha + (1 - \alpha) j/4$, for $j = 1, 2, 3$. The script returns values $(\mu, \sigma)$ which satisfy:
\[
Q_j = \exp{(\mu + \sigma x_j)}, \quad j = 1, 2, 3.
\]
The output ignores Poissonization. Script  \ref{sc:lognormalestim} does not solve our problem, because
$\alpha$ is unknown. Our approach in the next two scripts is to ``guess'' $\alpha$,
then see what quartiles would be generated from the resulting $(\mu, \sigma)$ estimate:
then evaluate the quality of our guess by matching these quartiles with the real ones.

\subsubsection{Mathematica script for \texttt{censoredQuartiles}}\label{sc:censoredQuartiles}
\begin{Verbatim}[frame=single,fontsize=\small]
censoredQuartiles[Mu_, Sigma_] := Module[{Alpha, y, z, Q1, Q2, Q3},
   y = RandomVariate[LogNormalDistribution[Mu, Sigma], 50000];
   z = Select[Map[RandomVariate[PoissonDistribution[#]] &, y], # > 0 &];
   Alpha = 1.0 - N[Length[z]/Length[y]]; 
   Q1 = InverseCDF[LogNormalDistribution[Mu, Sigma], 0.75*Alpha + 0.25];
   Q2 = InverseCDF[LogNormalDistribution[Mu, Sigma], 0.5*Alpha + 0.5];
   Q3 = InverseCDF[LogNormalDistribution[Mu, Sigma], 0.25*Alpha + 0.75];
   Return[{Q1, Q2, Q3}];
   ];
\end{Verbatim}
\subsubsection{Parsing the script \ref{sc:censoredQuartiles} }
Script \ref{sc:censoredQuartiles} is a partial inverse of script \ref{sc:lognormalestim}.
The arguments are $\mu$ and $\sigma$. A sample $\mathbf{y}$ of 50K Log-normal$(\mu, \sigma)$ random variates is
generated, and then 50K independent Poisson variates $\mathbf{z}$ where $z_i$ has mean $y_i$.
The proportion for which $z_i=0$ is assigned to $\alpha$.  The script return the quartiles of the
censored Log-normal$(\mu, \sigma)$ distribution.

\subsubsection{Mathematica script for \texttt{logRatioError}}\label{sc:logRatioError}
This script uses both scripts \ref{sc:lognormalestim} and \ref{sc:censoredQuartiles} to evaluate
how well our ``guess`'' for $\alpha$ gives realistic quartilies.
\begin{Verbatim}[frame=single,fontsize=\small]
logRatioError[q1t_, q2t_, q3t_,Alpha_] := 
  Module[{Mu, Sigma, q1p, q2p, q3p, logRatios },
   {Mu, Sigma} = lognormalParams[q1t, q2t, q3t, Alpha]; 
   {q1p, q2p, q3p} = censoredQuartiles[Mu, Sigma]; 
   logRatios = Log[N[{q1p, q2p, q3p} /{q1t, q2t, q3t}]]; 
   Return[Max[Abs[logRatios] ] ]; 
   ];
\end{Verbatim}
\subsubsection{Parsing the script \ref{sc:logRatioError}}
View the output of script \ref{sc:censoredQuartiles} as $(Q_1(\alpha), Q_2(\alpha), Q_3(\alpha))$ 
to show the dependence on $\alpha$. Script \ref{sc:logRatioError}
quantifies the deviation of these quartiles, shown as
\texttt{q1p, q2p, q3p}, from the observed ones $Q_1, Q_2, Q_3$, 
represented in the script as \texttt{q1t, q2t, q3t}, using the function
\begin{equation}\label{e:ratioerror}
    \max_{i = 1, 2, 3}  \left(\left|\log{ \frac{Q_i(\alpha)}{Q_i} }\right|\right),
\end{equation}

\subsubsection{Poisson Log-normal parameter fit by combining scripts \ref{sc:lognormalestim}, 
\ref{sc:censoredQuartiles}, \ref{sc:logRatioError} }
\begin{figure}
	\caption{\textbf{Fitting a censoring rate }
	\textit{In order to fit a Poisson Log-normal distribution when zero Poisson values are censored,
 we search for the censoring rate (horizontal axis) which minimizes deviation (\ref{e:ratioerror})
on vertical axis of implied quartiles from the observed ones.
		}
	} \label{f:ratioerror}	
	\begin{center}
		\scalebox{0.5}{\includegraphics{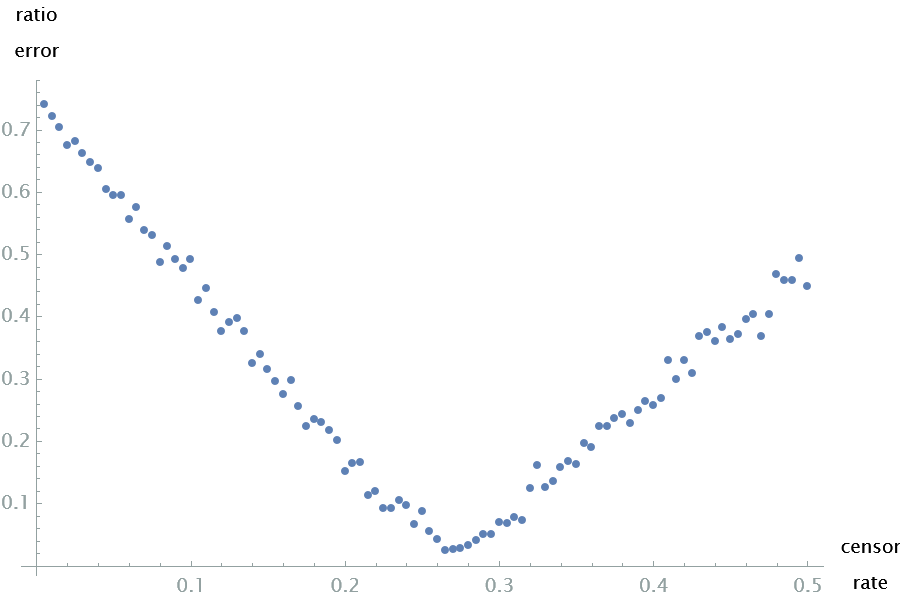}} 
	\end{center}	
\end{figure}

Choose a small $\delta>0$ and
iterate over values for $\alpha$ of the form $\delta, 2 \delta, 3 \delta, \ldots$, computing
a $(\mu, \sigma)$ pair for each with script 
\ref{sc:lognormalestim}. For each pair, script \ref{sc:censoredQuartiles}
returns quartiles based on a 50K Monte Carlo estimation of the rate of censoring.
Exhaustive search over multiples of $\delta$ reveals the value of $\alpha$
which minimizes (\ref{e:ratioerror}) using script \ref{sc:logRatioError}.
Call this optimal value $\alpha_*$.
The minimization is illustrated in Figure \ref{f:ratioerror}.
Finally apply the script \ref{sc:lognormalestim} with arguments $Q_1, Q_2, Q_3, \alpha_*$ to obtain a fitted pair $(\mu, \sigma)$
of Log-normal parameters which take account of censoring of zeros.

Taking $\delta:=0.005$ and $(Q_1, Q_2, Q_3):=(9, 109, 2559)$, this led to a fitted
value of $\alpha:=0.295$, and fitted Log-normal parameters $\mu =2.556$, $\sigma = 5.640$,
which are not far from the true values $(\mu, \sigma)=(2.6, 5.6)$.
The objective function is noisy because of sampling error in the simulation
occurring in the \texttt{censoredQuartiles} function, but nevertheless the minimum is sharply defined;
see Figure \ref{f:ratioerror}.
This concludes our quartile-based fitting of
the zero-censored Poisson Log-normal model.

In the case where
Poisson random variable $Z$ is thinned at rate $t \in [0,1]$, the conditional mean will become 
$t e^{\mu + \sigma Y} = e^{\mu' + \sigma Y}$ where $\mu' = \mu + \log{t}$. The thinning rate $t$ will be confounded with $\mu$,
but estimation of $\sigma$ is not affected.

\end{document}